\documentclass{pspum-l}
\usepackage[active]{srcltx}

\usepackage{graphicx}
\usepackage{hyperref}
\usepackage{enumerate}
\usepackage{bbm}

\newtheorem{thm}{Theorem}[section]
\newtheorem{lem}[thm]{Lemma}%
\newtheorem{prop}[thm]{Proposition}%
\theoremstyle{remark}
\theoremstyle{plain}

\numberwithin{equation}{section}

\def\QQ{{\mathbb Q}}
\def\PP{{\mathbb P}}
\def\RR{{\mathbb R}}

\def\ZZ{{\mathbb Z}}

\def\pp{{\mathbbm p}}

\def\veca{{\text{\boldmath$a$}}}
\def\vecb{{\text{\boldmath$b$}}}

\def\vece{{\text{\boldmath$e$}}}
\def\vech{{\text{\boldmath$h$}}}

\def\vecell{{\text{\boldmath$\ell$}}}

\def\vecn{{\text{\boldmath$n$}}}
\def\vecq{{\text{\boldmath$q$}}}
\def\vecQ{{\text{\boldmath$Q$}}}
\def\vecp{{\text{\boldmath$p$}}}
\def\vecP{{\text{\boldmath$P$}}}

\def\vecs{{\text{\boldmath$s$}}}

\def\vect{{\text{\boldmath$t$}}}

\def\vecv{{\text{\boldmath$v$}}}
\def\vecV{{\text{\boldmath$V$}}}
\def\vecw{{\text{\boldmath$w$}}}

\def\vecx{{\text{\boldmath$x$}}}

\def\vecy{{\text{\boldmath$y$}}}

\def\vecbeta{{\text{\boldmath$\beta$}}}

\def\veceta{{\text{\boldmath$\eta$}}}
\def\vecomega{{\text{\boldmath$\omega$}}}

\def\vecxi{{\text{\boldmath$\xi$}}}

\def\vecnull{{\text{\boldmath$0$}}}

\def\scrA{{\mathcal A}}
\def\scrB{{\mathcal B}}
\def\scrC{{\mathcal C}}
\def\scrD{{\mathcal D}}
\def\scrE{{\mathcal E}}
\def\scrF{{\mathcal F}}

\def\scrK{{\mathcal K}}
\def\scrL{{\mathcal L}}

\def\scrQ{{\mathcal Q}}
\def\scrP{{\mathcal P}}

\def\scrV{{\mathcal V}}
\def\scrW{{\mathcal W}}

\def\fO{{\mathfrak O}}

\def\fZ{{\mathfrak Z}}

\def\e{\mathrm{e}}

\def\id{\operatorname{id}}
\def\intl{{\operatorname{int}}}

\def\C{\operatorname{C{}}}

\def\L{\operatorname{L{}}}
\def\M{\operatorname{M{}}}

\def\S{\operatorname{S{}}}

\def\SL{\operatorname{SL}}
\def\ASL{\operatorname{ASL}}

\def\SO{\operatorname{SO}}

\def\T{\operatorname{T{}}}

\def\vol{\operatorname{vol}}

\def\GamG{\Gamma\backslash G}

\def\SLdZ{\SL(d,\ZZ)}
\def\SLdR{\SL(d,\RR)}

\def\trans{\,^\mathrm{t}\!}

\def\nbar{\overline{n}}
\def\xibar{\overline{\xi}}
\def\sigmabar{\overline{\sigma}}

\def\hatq{{\widehat{\vecq}}}
\def\hatp{{\widehat{\vecp}}}

\def\hatw{{\widehat{\vecw}}}

\def\hatQ{{\widehat{\vecQ}}}
\def\hatP{{\widehat{\vecP}}}

\title{Kinetic limits of dynamical systems}
\author{Jens Marklof}
\address{School of Mathematics, University of Bristol,
Bristol BS8 1TW, U.K.\newline
\rule[0ex]{0ex}{0ex} \hspace{8pt}{\tt j.marklof@bristol.ac.uk}}
\date{\today}

\begin{document}

\begin{abstract}
Since the pioneering work of Maxwell and Boltzmann in the 1860s and 1870s, a major challenge in mathematical physics has been the derivation of macroscopic evolution equations from the fundamental microscopic laws of classical or quantum mechanics. Macroscopic transport equations lie at the heart of many important physical theories, including fluid dynamics, condensed matter theory and nuclear physics. The rigorous derivation of macroscopic transport equations is thus not only a conceptual exercise that establishes their consistency with the fundamental laws of physics: the possibility of finding deviations and corrections to classical evolution equations makes this subject both intellectually exciting and relevant in practical applications. The plan of these lectures is to develop a renormalisation technique that will allow us to derive transport equations for the kinetic limits of two classes of simple dynamical systems, the Lorentz gas and kicked Hamiltonians (or linked twist maps). The technique uses the ergodic theory of flows on homogeneous spaces (homogeneous flows for short), and is based on joint work with Andreas Str\"ombergsson. 
\end{abstract}

\maketitle


\section{Motivation}\label{sec1}

It is perhaps surprising that, more than a century after Boltzmann's revolutionary ideas, we still don't have a complete understanding of the kinetic theory of the hard sphere gas. The only significant rigorous result to-date is due to Lanford \cite{Lanford1974}, who showed that, in the low-density limit {\em (Boltzmann-Grad limit)}, the evolution of the hard sphere gas converges to the solution of the Boltzmann equation for times shorter than the mean collision time; cf.~also the recent papers by Gallagher et al.~\cite{Gallagher:2012wq} and Pulvirenti et al.~\cite{Pulvirenti:2013un}, where complete derivations, variations and extensions are discussed. 
One of the difficulties of the hard sphere gas is of course that the number of degrees of freedom in a given volume tends to infinity as we approach the Boltzmann-Grad limit. Such difficulties are mirrored in the Boltzmann equation, whose solutions are still hard to analyse, despite the ground-breaking work of DiPerna and Lyons \cite{DiPerna:1989wl}.
The plan for these lectures is to study the Boltzmann-Grad limit of two classes of much simpler systems, whose limit kinetic equation is linear. Despite the simplicity of the setting, we will see that our analysis requires rather modern mathematical tools, and leads to some unexpected answers. In particular, the macroscopic transport equations that describe the dynamics in the Boltzmann-Grad limit may in general differ from the expected linear Boltzmann equation. This is due to subtle correlations between subsequent particle collisions. 
The systems we will focus on are:

(a) {\em The Lorentz gas}, a gas of non-interacting particles moving in a fixed array of spherical scatterers. Since the particles are non-interacting, the dynamics reduces to a one-particle dynamics. The Lorentz gas was first introduced by Lorentz \cite{Lorentz1905} in 1905 to model the motion of electrons in a crystal, and has since served as a fundamental mathematical model to study chaotic diffusion \cite{Bunimovich:1980ur,Bleher:1992ku,Szasz:2007uo,super}.

(b) {\em Kicked Hamiltonians.} The kicked rotor is the classic example, and the corresponding Chirikov standard map is one of the key models in chaos theory \cite{MR1169466}. Other examples of maps that are related to kicked Hamiltonians are the {\em linked twist maps} \cite{Sturman2006}.

\section{The Lorentz gas}\label{sec2}

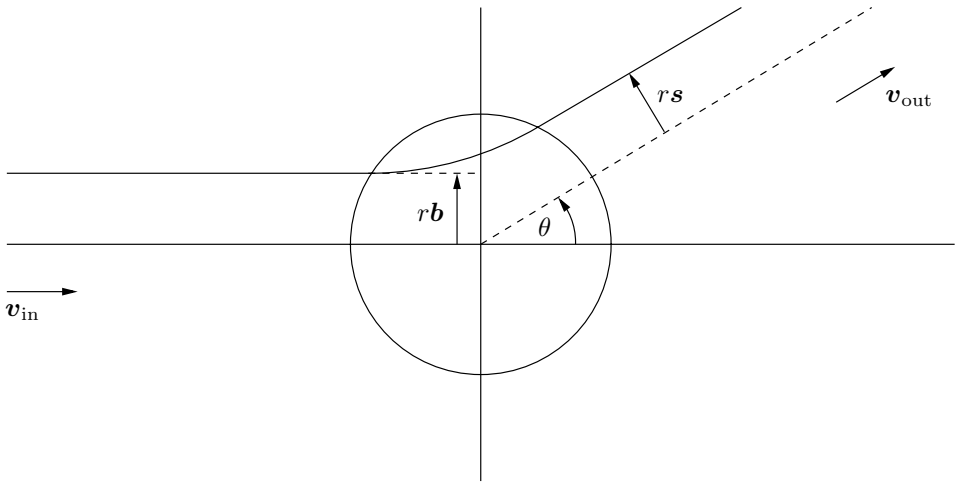
\begin{figure}
\begin{center}
\begin{minipage}{\textwidth}
\unitlength0.1\textwidth
\begin{picture}(10,5)(0,0)
\put(0,0){\includegraphics[width=\textwidth]{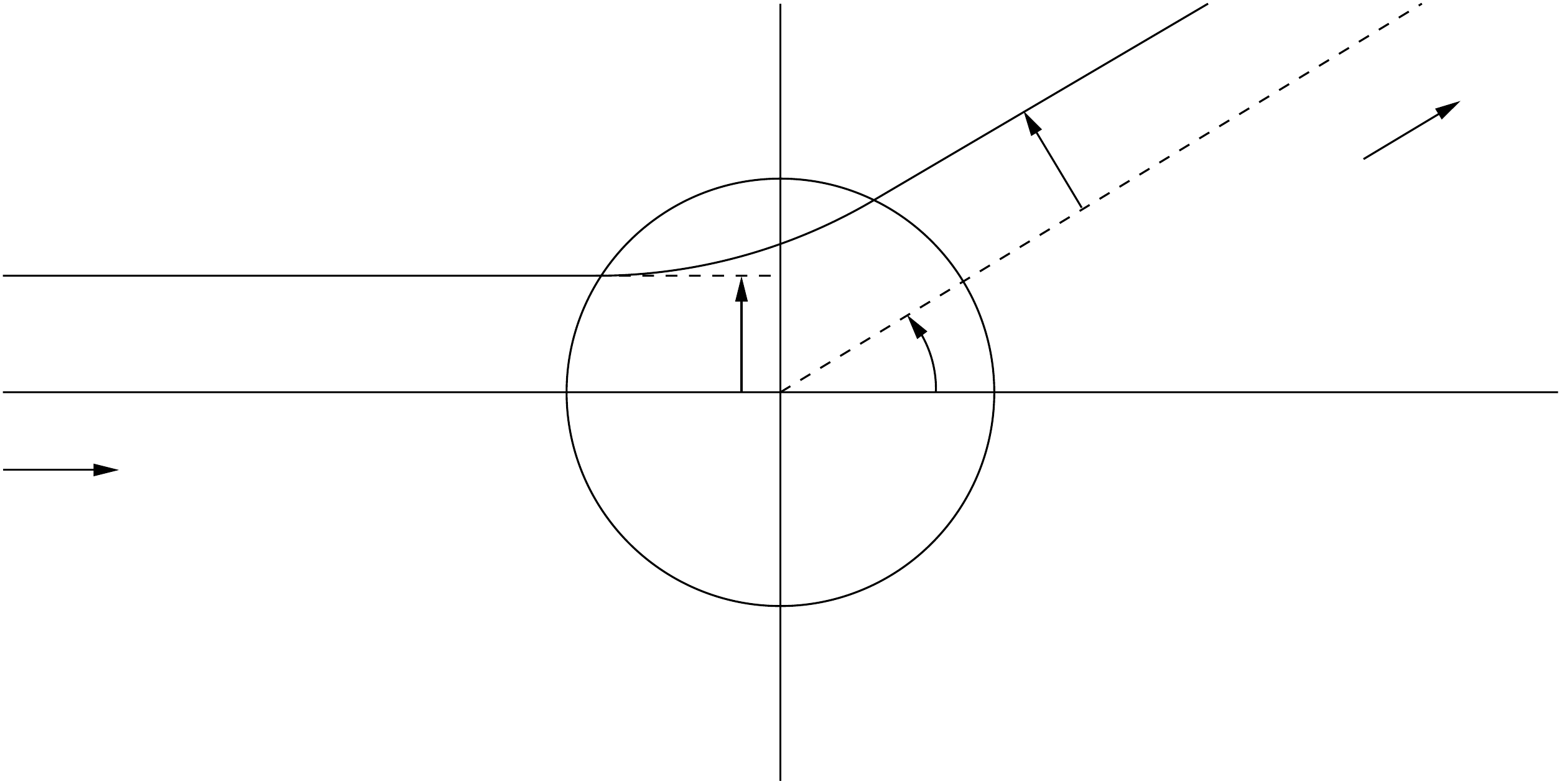}}
\put(0,1.75){$\vecv_\text{in}$} 
\put(4.33,2.75){$r\vecb$}
\put(5.6,2.6){$\theta$} 
\put(6.85,4){$r\vecs$}
\put(9.25,4){$\vecv_\text{out}$}
\end{picture}
\end{minipage}
\end{center}
\caption{Scattering in a sphere of radius $r$.} \label{fig1}
\end{figure}

Denote by $\scrP\subset\RR^d$ a point set. We assume throughout these lectures that $\scrP$ has constant density $\nbar>0$, i.e.\ for any $\scrD\subset\RR^d$ with $\vol(\partial\scrD)=0$ ($\vol$ denotes the Lebesgue measure in $\RR^d$ and $\partial\scrD$ the boundary of $\scrD$) we have
\begin{equation}\label{density}
\lim_{T\to\infty} \frac{\#(\scrP\cap T\scrD)}{\vol(T\scrD)} =\nbar.
\end{equation}
At each point in $\scrP$ we place a spherical scatterer of radius $r>0$, and consider the dynamics of a point particle in this infinite array of scatterers. We assume that the scatterers do not overlap. The scattering process at each scatterer is identical. It is furthermore assumed to be elastic (so that the particle speed after scattering is the same as before) and spherically symmetric (so that angular momentum is preserved). Since the particle speed is preserved, we may assume without loss of generality that between scattering events the speed is $\|\vecv\|=1$. We can then describe the scattering process in terms of a coordinate system where (cf.~Figure~\ref{fig1})
\begin{equation}
\vecv_\text{in} =\vece_1:= (1, 0 , \ldots , 0 ).
\end{equation}
(All vectors are represented as row vectors.)
The {\em impact parameter} $\vecb$ is the orthogonal projection of the point of impact onto the plane orthogonal to $\vecv_\text{in}$, measured in units of $r$. In the present coordinate system, we have $\vecb=(0,\vecw)$ for some $\vecw\in\scrB_1^{d-1}$ (we will also refer to $\vecw$ as impact parameter). The outgoing velocity is for $\vecb\neq\vecnull$
\begin{equation}\label{inout}
\vecv_\text{out} =  \vecv_\text{in} \cos\theta +  (0,\hatw) \sin\theta  ,
\end{equation}
where the angle $\theta$ is called the {\em scattering angle} and $\hatw:=w^{-1}\vecw$ with $w:=\|\vecw\|$. For $\vecw=\vecnull$ we simply assume $\vecv_\text{out} =  -\vecv_\text{in}$. By the assumed spherical symmetry, $\theta=\theta(w)$ is only a function of the length $w$ of the impact parameter $\vecw$. Equation \eqref{inout} can be expressed as
\begin{equation}\label{inout2}
\vecv_\text{out} =  \vecv_\text{in} S(\vecw)^{-1}  ,
\end{equation}
with the matrix
\begin{equation}\label{SbLG}
S(\vecw) = \exp\begin{pmatrix} 0 & -\theta(w) \hatw \\ \theta(w) \trans\hatw & 0_{d-1} \end{pmatrix}\in\SO(d).
\end{equation}
The {\em exit parameter} is defined as the orthogonal projection the point of exit onto the plane orthogonal to $\vecv_\text{out}$, and is given by 
\begin{equation}\label{inout2exit}
\vecs = -w \vecv_\text{in} \sin\theta +  (0,\vecw) \cos\theta = (0,\vecw) S(\vecw)^{-1}  .
\end{equation}

We assume that one of the following condition holds:

\begin{itemize}
\item[(A)] $\theta\in\C^1([0,1))$ is strictly decreasing with $\theta(0)=\pi$ and $\theta(w)> 0$. 
\item[(B)] $\theta\in\C^1([0,1))$ is strictly increasing with $\theta(0)=-\pi$  and $\theta(w)< 0$. 
\end{itemize}

This hypothesis is satisfied for many scattering maps, e.g.~specular reflection (where $\theta(w)=\pi-2\arcsin(w)$ so (A) holds) or the scattering in the muffin-tin Coulomb potential $V(\vecq)=\alpha \max(r \|\vecq\|^{-1}-1,0)$ with $\alpha\notin \{-2E,0\}$, where $E$ denotes the
total energy, cf.~\cite{Marklof:2011ho}. The above assumption implies in particular that, for any given $\vecv_\text{in}$, the map $\vecb\mapsto \vecv_\text{out}$ is invertible. We denote by $\sigma(\vecv_\text{in},\vecv_\text{out})$ the Jacobian of the inverse map:
\begin{equation}\label{inout777}
\sigma(\vecv_\text{in},\vecv_\text{out}) \,d\vecv_\text{out} = d\vecb,
\end{equation}
where $d\vecb$ is the Lebesgue measure on the hyperplane orthogonal to $\vecv_\text{in}$.
This Jabobian is called the {\em differential cross section}. The {\em total scattering cross section} in these units is \begin{equation}
\sigmabar:=\int_{\S_1^{d-1}} \sigma(\vecv_\text{in},\vecv_\text{out}) \,d\vecv_\text{out} =\vol\scrB_1^{d-1}.
\end{equation}

\begin{figure}
\begin{center}
\begin{minipage}{0.5\textwidth}
\unitlength0.1\textwidth
\begin{picture}(10,9)(0,0)
\put(0,0){\includegraphics[width=\textwidth]{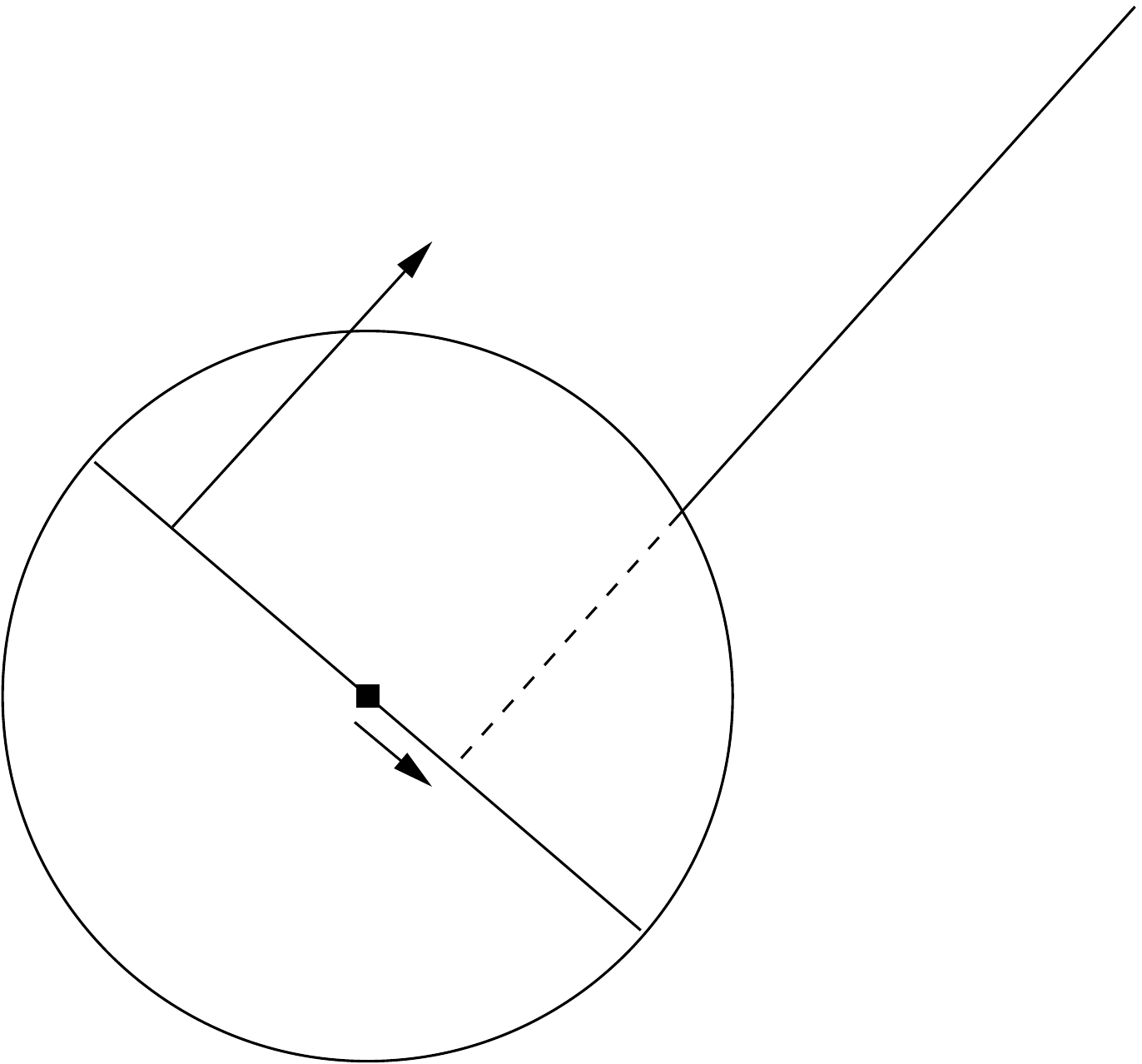}}
\put(2.7,2.2){$\vech$} 
\put(3.5,3.5){$\vecy$}
\put(2.7,5.5){$\vecv$} 
\put(8,6.5){$\ell$}
\end{picture}
\end{minipage}
\end{center}
\caption{Coordinates relative to a scatterer.} \label{fig2}
\end{figure}
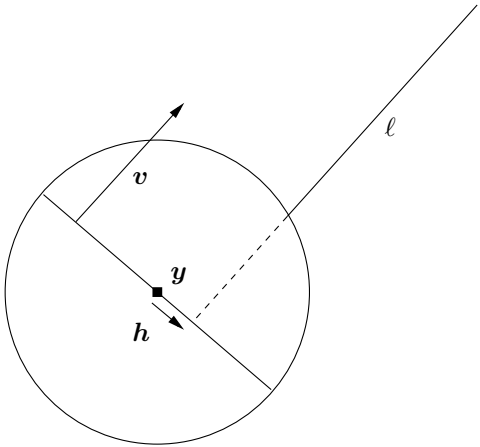

Outside the scatterers the particle moves along straight lines with constant velocity $\vecv$ and $\|\vecv\|=1$. It will be convenient to ignore the dynamics inside the scatterer and assume the scattering is instantaneous. The configuration space of the Lorentz gas is then
\begin{equation}
\scrK_r= \RR^d \setminus (\scrP+\scrB_r^d)
\end{equation}
where $\scrB_r^d$ is the ball in $\RR^d$ of radius $r$. Its phase space is $\T^1(\scrK_r)$, the unit tangent bundle of $\scrK_r$. We use the convention that for $\vecq\in\partial\scrK_r$ the vector $\vecv$ points away from the scatterer, so that $\vecv$ describes the velocity {\em after} the collision. The Liouville measure for the dynamics is $d\vecq\,d\vecv$, where $\vecq\in\scrK_r$ denotes the particle position.
We introduce an alternative parametrisation of phase space by setting (see Figure \ref{fig2})
\begin{equation}
(\vecq,\vecv)=(\vech+\ell\vecv,\vecv) .
\end{equation}
Here $\vech\in\scrB_{r,\vecv}^{d-1}+\vecy$ where $\scrB_{r,\vecv}^{d-1}$ is the hyperdisk of radius $r$ perpendicular to $\vecv$, and $\vecy\in\scrP$ is the scatterer location. The quantity $\ell\geq \sqrt{r^2-\|\vech\|^2}$ describes the distance travelled since the leaving the last scatterer. In these coordinates we have (for arbitrary fixed $\vecv$)
\begin{equation}
d\vecq = d\vech\,d\ell .
\end{equation}

\section{Mean free path length in the Lorentz gas}\label{sec3}

The mean free path length is the ``average'' distance travelled between collisions. As we are travelling with speed $\|\vecv\|=1$, the mean free path length is the same as the mean collision time. For a particle with initial data $(\vecq,\vecv)\in\T^1(\scrK_r)$, denote by 
\begin{equation}
\tau_1(\vecq,\vecv)=\inf\{ t>0 : \vecq+t\vecv\notin\scrK_r\}
\end{equation}
the first collision time. The collision time for a particle which previously hit the scatterer at $\vecy\in\scrP$ is then, averaged with respect to $\vech$ and $\vecv$,
\begin{equation}
\frac{1}{ \vol\S_1^{d-1}}\;\frac{1}{\vol\scrB_r^{d-1}} \int_{\S_1^{d-1}}\int_{\scrB_{r,\vecv}^{d-1}}
\tau_1(\vech+\vecy+\vecv \sqrt{r^2-\|\vech\|^2}, \vecv)\, d\vech\,d\vecv .
\end{equation}
In order to define the mean collision time/path length for the full Lorentz gas with scatterer configuration $\scrP$, we fix a convex bounded set $\scrD\subset\RR^d$, and consider only scatterers inside the set $T\scrD$ with $T$ large. The scatterer locations are given by the set
\begin{equation}
\scrP_T := \{ \vecy\in\scrP: \vecy+\scrB_r^d\subset T\scrD \},
\end{equation}
and in view of \eqref{density} we have $\#\scrP_T\sim \nbar \vol(T\scrD)$. We now define the mean free path length for the finite Lorentz gas in $T\scrD$, where we treat the boundary of $T\scrD$ in the same way as a scatterer (the first collision time refers now to a collision with a scatterer or the boundary $\partial(T\scrD)$), and take the limit $T\to\infty$. If the answer is independent of $\scrD$, we have a well defined mean free path length for the infinite Lorentz gas. Truncations of this kind are a standard trick in statistical mechanics.

Let us denote by $\scrD_\vecv$ the orthogonal projection of $\scrD$ onto the hyperplane perpendicular to $\vecv$. The mean collision time is now
\begin{equation}\label{mean1}
\begin{split}
& \frac{1}{ \vol\S_1^{d-1}} \int_{\S_1^{d-1}} \frac{1}{\#\scrP_T \vol\scrB_r^{d-1}+ T^{d-1}\vol\scrD_\vecv} 
\times \\ 
& \times
\bigg\{ \sum_{\vecy\in\scrP_T}  \int_{\scrB_{r,\vecv}^{d-1}}
\tau_1(\vech+\vecy+\vecv \sqrt{r^2-\|\vech\|^2}, \vecv)\, d\vech \\
& + \int_{T\scrD_\vecv} \tau_1(\vecq(\vech),\vecv) \, d\vech \bigg\}  d\vecv ,
\end{split}
\end{equation}
where, for given $\vecv$ and $T$, the function $T\scrD_\vecv\to\partial (T\scrD)$, $\vech\mapsto\vecq(\vech)$, parametrises the points of $\vecq\in\partial (T\scrD)$ whose inward pointing normal $\vecn$ satisfies $\angle(\vecv,\vecn)\leq\pi/2$. We write
\begin{equation}
\tau_1=\int_0^{\tau_1} d\ell
\end{equation}
and recall $d\vecq = d\vech\,d\ell$ for $\vecv$ fixed. Thus, \eqref{mean1} equals
\begin{equation}
 \frac{1}{ \vol\S_1^{d-1}} \int_{\S_1^{d-1}}  \frac{\vol(T\scrD)-\#\scrP_T\vol\scrB_r^{d}}{\#\scrP_T\vol\scrB_r^{d-1}+ T^{d-1}\vol\scrD_\vecv} \, d\vecv,
\end{equation}
which is asymptotic to (as $T\to\infty$)
\begin{equation}
\frac{\vol(T\scrD)-\#\scrP_T\vol\scrB_r^d}{\#\scrP_T \vol\scrB_r^{d-1}}
\longrightarrow 
\frac{1-\nbar \vol\scrB_r^d}{\nbar \vol\scrB_r^{d-1}}
=\frac{1-\nbar r^d \vol\scrB_1^d}{\nbar r^{d-1} \vol\scrB_1^{d-1}}.
\end{equation}
Note that this expression for the mean free path length of the Lorentz gas is indeed independent of the choice of $\scrD$. It is also independent of the scatterer configuration $\scrP$, given its density is $\nbar$.

\section{The Boltzmann-Grad limit of the Lorentz gas}\label{sec4}

In the case of the Lorentz gas, the Boltzmann-Grad limit is defined as the limit of low scatterer density (as opposed to the limit of low particle density in Boltzmann's hard sphere gas). {\em Density} refers here to the {\em volume} density, i.e., the relative volume $\nbar r^d \vol\scrB_1^d$ occupied by the scatterers, rather than their {\em number} density $\nbar$. For a fixed scatter configuration $\scrP$ (and thus fixed $\nbar$) the Boltzmann-Grad limit corresponds therefore to taking $r\to 0$. From Section \ref{sec3} we infer that in this case the mean free path length is asymptotically
\begin{equation}
\frac{1}{\nbar r^{d-1} \vol\scrB_1^{d-1}} .
\end{equation}
To capture the dynamics of the Lorentz gas in the Boltzmann-Grad limit $r\to 0$, we measure length and time in units of the mean free path length/mean collision time. This is achieved by using the {\em macroscopic} coordinates 
\begin{equation}
(\vecQ(t),\vecV(t)) = (r^{d-1} \vecq(r^{-(d-1)} t), \vecv(r^{-(d-1)} t) ) .
\end{equation}

Denote the corresponding macroscopic configuration space by $\widetilde\scrK_r=r^{d-1}\scrK_r$. 
In these units, the mean free path length and mean collision time are equal to 
\begin{equation}\label{MFP}
\xibar=\frac{1}{\nbar\,\sigmabar} ,
\end{equation}
with the total scattering cross section $\sigmabar = \vol\scrB_1^{d-1}$.

The evolution of an initial macroscopic particle density $f\in\T^1(\widetilde\scrK_r)$ is defined by the linear operator 
\begin{equation}
[L_r^t f](\vecQ,\vecV) := f(\vecQ(-t),\vecV(-t))
\end{equation}
where $(\vecQ(t),\vecV(t))$ is the solution of Hamilton's equations with initial condition $(\vecQ(0),\vecV(0))=(\vecQ,\vecV)$. We extend $L_r^t$ to a linear operator on $\T^1(\RR^d)$ by setting 
\begin{equation}
[L_r^t f](\vecQ,\vecV) := f(\vecQ,\vecV) \qquad \text{if $\vecQ\notin\widetilde\scrK_r$.}
\end{equation}
We would like to answer the following questions:

(1) For a given scatterer configuration $\scrP$, does $L_r^t$ have a (weak) limit as $r\to 0$? That is, for every $t>0$ is there 
\begin{equation}
L^t:\L^1(\T^1(\RR^d))\to\L^1(\T^1(\RR^d))
\end{equation}
such that for every ``nice'' bounded $\scrA\subset\T^1(\RR^d)$ and $f\in\L^1(\T^1(\RR^d))$
\begin{equation}
\lim_{r\to 0} \int_\scrA L_r^t f(\vecQ,\vecV) \, d\vecQ\,d\vecV = \int_\scrA L^t f(\vecQ,\vecV) \, d\vecQ\,d\vecV \; ?
\end{equation}

(2) More generally, is there a random flight process describing the dynamics in the Boltzmann-Grad limit and what are its properties?

In 1905 Lorentz \cite{Lorentz1905} answered these questions using Boltzmann's heuristics. They were later confirmed rigorously for random scatterer configurations $\scrP$ (such as a typical realisation of a Poisson process) by Gallavotti \cite{Gallavotti1969}, Spohn \cite{Spohn:1978tt} and Boldrighini, Bunimovich and Sinai \cite{Boldrighini:1983jm}. In this case $L^t$ exists and $f_t:=L^t f$ satisfies the {\em linear} Boltzmann equation (also referred to as the {\em kinetic Lorentz equation})
\begin{equation}\label{LBeq}
(\partial_t +\vecV\cdot\partial_\vecQ) f_t(\vecQ,\vecV) 
= \nbar \int_{\S_1^{d-1}} [f_t(\vecQ,\vecV')-f_t(\vecQ,\vecV)] \,\sigma(\vecV,\vecV')\, d\vecV'.
\end{equation}
As we shall see, the linear Boltzmann equation fails for other scatterer configurations, and the Boltzmann-Grad limit leads to a more complicated limit process. This process is governed by a transport equation which, aside from position $\vecQ$ and velocity $\vecV$, also requires the following data:
\begin{itemize}
\item $\xi$ is the distance to the next collision (in the above macroscopic coordinates)
\item $\vecomega$ are variables that characterise the next scattering event. One example is the impact parameter $\vecb$, but for some scatterer configurations $\scrP$ more information will be required (as we shall see below). We denote by $\pp$ the probability measure on the relevant parameter space $\Omega$. 
\end{itemize}
The transport equation reads
\begin{equation}\label{FPK}
(\partial_t +\vecV\cdot\partial_\vecQ-\partial_\xi) f_t(\vecQ,\vecV,\xi,\vecomega) \\
= [\scrC f_t](\vecQ,\vecV,\xi,\vecomega)
\end{equation}
with the collision operator
\begin{equation}
 [\scrC f](\vecQ,\vecV,\xi,\vecomega) = \int_\Omega k(\vecomega_+',\xi,\vecomega) f(\vecQ,\vecV'(\vecomega_+',\vecV),0,\vecomega'(\vecomega_+',\vecV))\,d\pp(\vecomega_+')
\end{equation}
where $\vecomega_+'$ denotes the hidden variable $\vecomega'$ with the impact parameter $\vecb'$ replaced by the exit parameter $\vecs'$ of the same collision, and
$\vecV'(\vecomega_+',\vecV)$ is the incoming velocity $\vecV'$ so that the outgoing velocity with exit parameter $\vecs'$ is $\vecV$; hence $\vecV'(\vecomega',\vecV)$ is only a function of $\vecs'$ and $\vecV$. Note that $\vecs'$ and the impact parameter $\vecb$ of the next collision are in the same hyperdisk orthogonal to $\vecV$.

The random flight process governed by equation \eqref{FPK} is that of a particle travelling with constant velocity $\vecV_{n-1}$ and a ``label'' $\vecomega_{n-1}$, which at random time $\xi_n$ after the last scattering changes to $\vecomega_n$ with probability $k(\vecomega_{n-1},\xi_n,\vecomega_n)\,d\xi_n\,d\pp(\vecomega_n)$. The particle changes direction, and its new velocity $\vecV_{n}=\vecV(\vecomega_n,\vecV_{n-1})$ is determined by the scattering map. The particle continues travelling with constant velocity $\vecV_{n}$ and label $\vecomega_{n}$, and again after time $\xi_{n+1}$ changes to $\vecomega_{n+1}$ with probability $k(\vecomega_{n},\xi_{n+1},\vecomega_{n+1})\,d\xi_{n+1}\,d\pp(\vecomega_{n+1})$, and so on. This random flight process is a continuous-time Markov process, and eq.~\eqref{FPK} is the corresponding {\em Fokker-Planck-Kolmogorov equation}. 

The time-reversibility of the underlying microscopic dynamics (for every fixed $r>0$) implies that the transition kernel $k$ is symmetric, i.e.
\begin{equation}\label{ksymm}
k(\vecomega,\xi,\vecomega') = k(\vecomega',\xi,\vecomega).
\end{equation}
If $\tilde f_0(\vecQ,\vecV)$ is the initial particle density, the ``physical'' initial condition of our transport equation is 
\begin{equation}
f_0(\vecQ,\vecV,\xi,\vecomega) = \tilde f_0(\vecQ,\vecV) K(\xi,\vecomega)
\end{equation}
where 
\begin{equation}\label{K}
K(\xi,\vecomega) := \frac{1}{\,\xibar\,} \int_\xi^\infty \int_\Omega k(\vecomega',\xi',\vecomega)\, d\xi'\, d\pp(\vecomega').
\end{equation}

Note that the choice
\begin{equation}
f_t(\vecQ,\vecV,\xi,\vecomega) = K(\xi,\vecomega)
\end{equation}
is a stationary solution of the transport equation, since
\begin{equation}\label{stats}
-\partial_\xi \int_\xi^\infty \int_\Omega k(\vecomega',\xi',\vecomega)\, d\xi'\, d\pp(\vecomega') =\xibar
\int_\Omega k(\vecomega',\xi,\vecomega)\, K(0,\vecomega') \, d\pp(\vecomega') .
\end{equation}
To prove \eqref{stats}, observe that the left hand side equals 
\begin{equation} \label{lhs}
\int_\Omega k(\vecomega',\xi,\vecomega)\, d\pp(\vecomega') .
\end{equation}
As to the right hand side,
\begin{equation}
\begin{split}
K(0,\vecomega') & = \xibar^{-1} \int_0^\infty \int_\Omega k(\vecomega'',\xi,\vecomega')\, d\xi\, d\pp(\vecomega'') \\
& = \xibar^{-1}  \int_0^\infty \int_\Omega k(\vecomega',\xi,\vecomega'')\, d\xi\, d\pp(\vecomega'') \\
& = \xibar^{-1} ,
\end{split}
\end{equation}
which shows that the right hand side of \eqref{stats} equals \eqref{lhs}.

We obtain the evolved particle density $\tilde f_t(\vecQ,\vecV)$ with initial data $\tilde f_0(\vecQ,\vecV)$ from the solution of \eqref{FPK} by 
\begin{equation}
\tilde f_t(\vecQ,\vecV) = \int_0^\infty \int_\Omega f_t(\vecQ,\vecV,\xi,\vecomega) \, d\xi\,d\pp(\vecomega) .
\end{equation}

The kernel $k(\vecomega',\xi,\vecomega)$ depends in general on the scatterer configuration $\scrP$, and in these lectures we will explore under which assumptions we can expect the Boltzmann-Grad limit to exist and to be described by the above kinetic equation \eqref{FPK}. 

Let us remark that 
\begin{equation}\label{FreePL}
\overline\Phi_0(\xi) = \int_\Omega\int_\Omega k(\vecomega',\xi,\vecomega) \, d\pp(\vecomega')\, d\pp(\vecomega).
\end{equation}
describes the distribution of free path length between consecutive collisions in the Boltzmann-Grad limit. We have
\begin{equation}
\int_0^\infty \overline\Phi_0(\xi)\, d\xi = 1.
\end{equation}
The mean free path length \eqref{MFP} satisfies the relation 
\begin{equation}
\begin{split}
\xibar & = \int_0^\infty \xi \overline\Phi_0(\xi)\, d\xi \\
& = \int_\Omega \int_\Omega\int_0^\infty  \xi k(\vecomega',\xi,\vecomega)\, d\xi\, d\pp(\vecomega')\, d\pp(\vecomega) \\
& = - \xibar \int_\Omega \xi K(\xi,\vecomega) \bigg|_{\xi=0}^\infty d\pp(\vecomega) + \xibar \int_\Omega \int_0^\infty K(\xi,\vecomega)\,d\xi\,d\pp(\vecomega) \\
& = \xibar \int_\Omega \int_0^\infty K(\xi,\vecomega)\,d\xi\,d\pp(\vecomega) ,
\end{split}
\end{equation}
which shows that $K(\xi,\vecomega)$ is a probability density with respect to $d\xi\,d\pp(\vecomega)$.

As a first example, assume the transition kernel is
\begin{equation}\label{ex0}
k(\vecomega',\xi,\vecomega) = \xibar^{-1} \e^{-\xi/\xibar} 
\end{equation}
and the hidden variables comprise only the impact parameter normalised so that $d\pp(\vecomega)$ is a probability measure:
\begin{equation}
\vecomega:=\vecb,\qquad d\pp(\vecomega):=\sigmabar^{-1} d\vecb ,
\end{equation}
where $d\vecb$ is the Lebesgue measure on the hyperplane orthogonal to $\vecV$.
Then \eqref{K} yields
\begin{equation}
K(\xi,\vecomega) = \xibar^{-1}  \e^{-\xi/\xibar} 
\end{equation}
and the ansatz
\begin{equation}
f_t(\vecQ,\vecV,\xi,\vecomega)=\tilde f_t(\vecQ,\vecV)\, K(\xi,\vecomega) 
\end{equation}
in the transport equation \eqref{FPK} yields
\begin{equation}
\begin{split}
& (\partial_t +\vecV\cdot\partial_\vecQ+\xibar^{-1})  \tilde f_t(\vecQ,\vecV) \\
& = \xibar^{-1} \int_{\scrB_1^{d-1}} \tilde f_t(\vecQ,\vecV'(\vecs',\vecV))\, \sigmabar^{-1} d\vecs'\\
& = \nbar \int_{\scrB_1^{d-1}} \tilde f_t(\vecQ,\vecV') \sigma(\vecV,\vecV')\, d\vecV', 
\end{split}
\end{equation}
which is the linear Boltzmann equation \eqref{LBeq}. Here we have used $d\vecs'=\sigma(\vecV,\vecV')\, d\vecV'$, which follows from the time-reversal of relation \eqref{inout777}, $d\vecs'=\sigma(-\vecV,-\vecV')\, d\vecV'$ and the spherical invariance of the scattering map,
\begin{equation}
\sigma(-\vecV,-\vecV')=\sigma(\vecV,\vecV') .
\end{equation}

We will see in Section \ref{sec11} that the above scenario corresponds precisely to a random, Poisson distributed scatterer configuration.

\section{Kicked Hamiltonians}\label{sec5}

The machinery we will develop for the Lorentz gas can also be applied to a different class of systems, {\em kicked Hamiltonians.} An important example is the kicked rotor whose dynamics is described by Chirikov's {\em standard map} \cite{MR1169466}. In the mathematics literature, kicked Hamiltonians also appear in the guise of {\em linked twist maps}, whose ergodic properties were studied mainly in the 1980s, see \cite{Sturman2006} and references therein.

Given a bi-infinite sequence of kicking times,
\begin{equation}
-\infty \leftarrow \ldots < t_{-2} < t_{-1} < t_0 < t_1 < t_2 < \ldots \rightarrow \infty 
\end{equation}
we consider the kicked Hamiltonian
\begin{equation}
H(\vecq,\vecp,t)=\frac{\|\vecp\|^2}{2} +V(\vecq) \sum_{m\in\ZZ} \delta(t-t_m)
\end{equation}
with $\vecp,\vecq\in\RR^n$. The potential $V$ is of the form
\begin{equation}
V(\vecq)= r \sum_{j\in\ZZ^n} W\left(\frac{\vecq-\vecx_j}{r}\right)
\end{equation}
with a sequence of scatterer locations $\vecx_j$ and $W$ differentiable and supported in some compact $\Sigma\subset\RR^n$. Hamilton's equations yield 
\begin{equation}
\dot\vecp =- \partial_\vecq H  = -\partial V(\vecq) \sum_{m\in\ZZ} \delta(t-t_m),
\end{equation}
\begin{equation}
\dot\vecq =  \partial_\vecp H  =\vecp .
\end{equation}
The solution $\vecp(t)$ will be discontinuous at the kicking times $t=t_m$; here we define
\begin{equation}
\vecp(t):=\lim_{\epsilon\to 0_+} \vecp(t+\epsilon)
\end{equation}
so $\vecp(t)$ is continuous from the right and represents the momentum just after the kick. 
With this, the solution of Hamilton's equations for our kicked Hamiltonian is for the initial condition $(\vecq_0,\vecp_0)$, with $t_0\leq 0 <t_1$ and $t_m\leq t < t_{m+1}$, 
\begin{equation}
(\vecq(t),\vecp(t))=\Phi^t(\vecq_0,\vecp_0)
\end{equation}
with 
\begin{equation}
\Phi^t = \Phi_0^{t-t_m}\circ S_V \circ\Phi^{t_m-t_{m-1}}\circ S_V \circ\ldots \circ S_V \circ \Phi^{t_1}
\end{equation}
and
\begin{equation}
\Phi_0^t (\vecq,\vecp) := (\vecq+ t\vecp,\vecp) ,
\end{equation}
\begin{equation}
\begin{split}
S_V (\vecq,\vecp) & := (\vecq,\vecp-\partial V(\vecq)) \\
& =
\begin{cases}
\left(\vecq,\vecp-\partial W\left(\frac{\vecq-\vecx_j}{r}\right)\right) & \text{if $\vecq\in\vecx_j+r\Sigma$,}\\
 (\vecq,\vecp) & \text{otherwise.}
\end{cases}
\end{split}
\end{equation}
The variable $\vecw=\frac{\vecq-\vecx_j}{r}\in\Sigma$ plays the role of the impact parameter. The momenta before ($\vecp_\text{in}$) and after the kick $(\vecp_\text{out})$ are related by
\begin{equation}
\vecp_\text{out} =\vecp_\text{in} -\partial W(\vecw) .
\end{equation}
We assume the map 
\begin{equation}\label{KHmap}
\begin{split}
\Sigma & \to -\partial W(\Sigma) \\
\vecw & \mapsto \vecp=-\partial W(\vecw)
\end{split}
\end{equation}
is invertible and define the ``scattering'' cross section $\sigma(\vecp_\text{out}-\vecp_\text{in})$ by
\begin{equation}\label{sigKH}
\sigma(\vecp)\, d\vecp=d\vecw.
\end{equation}
The total scattering cross section is $\sigmabar=\vol\Sigma$.

\section{Geometric representation}\label{sec6}

\begin{figure}
\begin{center}
\begin{minipage}{0.7\textwidth}
\unitlength0.1\textwidth
\begin{picture}(10,8)(0,0)
\put(0,0){\includegraphics[width=\textwidth]{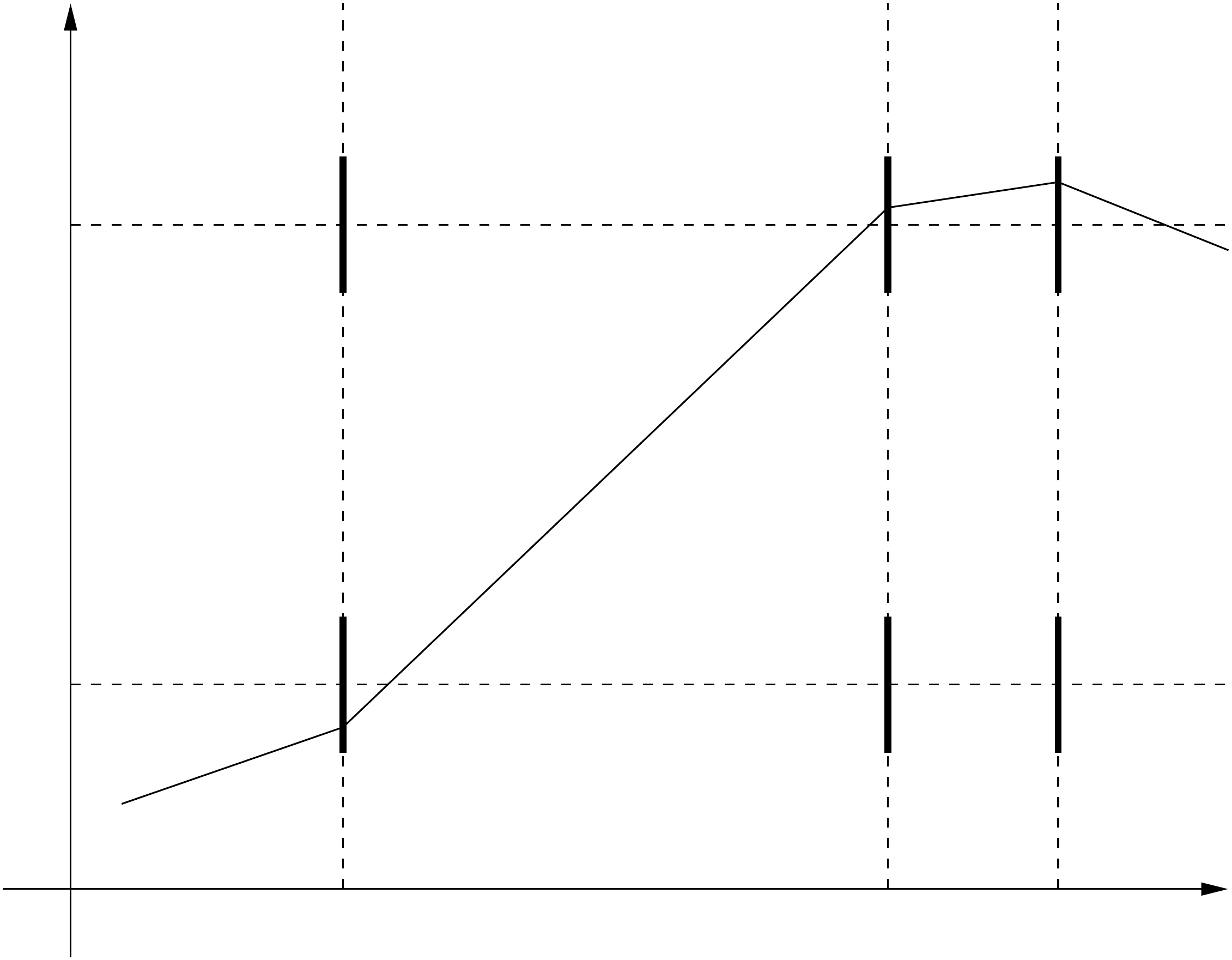}}
\put(0,2.2){$\vecx_1$} 
\put(0,5.85){$\vecx_2$} 
\put(0,7.5){$\vecq$} 
\put(2.7,0){$t_1$}
\put(7.1,0){$t_2$} 
\put(8.5,0){$t_3$}
\put(10,0){$q_0$}
\end{picture}
\end{minipage}
\end{center}
\caption{Trajectory in the geometric representation of a kicked Hamiltonian.} \label{fig3}
\end{figure}
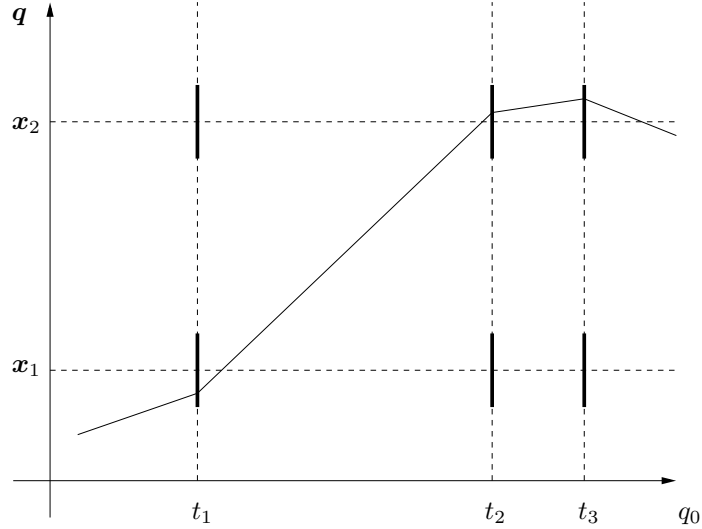

The following construction will show the close relationship between kicked Hamiltonians and the Lorentz gas.  Given the momentum $\vecp\in\RR^n$ and position $\vecq\in\RR^n$, we set
\begin{equation}
\hatq=(q_0,\vecq), \qquad \hatp=(p_0,\vecp)\in\RR^d \qquad \text{with $d=n+1$,}
\end{equation}
and define the set of scatterer locations in $\RR^d$ as
\begin{equation}\label{specialP}
\scrP=\{ (t_m,\vecx_j)  : m\in\ZZ,\; j\in\ZZ^n \} ,
\end{equation}
with $t_m,\vecx_j$ as in the previous Section \ref{sec5}.
Instead of spherical scatterers we now use (cf.~Figure \ref{fig3})
\begin{equation}
\{0\} \times r\Sigma .
\end{equation}
The dynamics between scattering events is (as in the Lorentz gas)
\begin{equation}
(\hatq(t),\hatp(t))=(\hatq(t')+(t-t')\hatp(t'),\hatp(t')),
\end{equation}
with $t_m \leq t'\leq t < t_{m+1}$. The scattering map is
\begin{equation}
\hatp_\text{in} = (p_0,\vecp_\text{in}) \mapsto \hatp_\text{out} = (p_0,\vecp_\text{out}) ,
\end{equation}
with $\vecp_\text{in}$, $\vecp_\text{out}$ as in Section \ref{sec5}, so that $p_0$ is preserved. 

Thus, instead of $\|\vecv\|$ as in the Lorentz gas, the constant of motion in this setting is $p_0$. For $p_0=1$ we recover precisely the dynamics of the kicked Hamiltonian. We will assume $p_0=1$ from now on. The phase space of the dynamics is thus $\RR^d\times\RR^{d-1}$, parametrised by $(\hatq,\vecp)$. The Liouville measure is $d\hatq\, d\vecp$ where $d\hatq=dq_0\, d\vecq$. The analogue of relation \eqref{inout2} is
\begin{equation}\label{inout222}
\hatp_\text{out} =  \hatp_\text{in} S(\vecw)^{-1}  ,
\end{equation}
where
\begin{equation}\label{SbKH}
S(\vecw)=\begin{pmatrix} 1 & \partial W(\vecw) \\ \trans\vecnull & 1_{d-1} \end{pmatrix}.
\end{equation}

There is nothing that prevents us from considering a more general point set than \eqref{specialP} (this allows for a different choice of potential locations $\vecx_j$ at each kicking time $t_m$). As in the Lorentz gas, we will for now only assume that the asymptotic density of $\scrP$ is $\nbar$, recall \eqref{density}. The configuration space is in this setting
\begin{equation}
\scrK_r=\RR^d\setminus (\scrP+\{0\}\times r\Sigma).
\end{equation}

\section{Mean collision time in kicked Hamiltonians}\label{sec7}

As in Section \ref{sec3}, we define the first collision time for initial condition $(\hatq,\hatp)$ by
\begin{equation}
\tau_1(\hatq,\hatp)=\inf\{ t>0 : \hatq+ t\hatp\notin\scrK_r\}.
\end{equation}
As in Section \ref{sec2}, we parametrise the phase space in terms of coordinates relative to the scatterers and the time travelled:
\begin{equation}
(\hatq,\hatp)= (\vech+\ell\hatp,\hatp)
\end{equation}
with $\vech\in\{0\}\times r\Sigma+\vecy$, $\vecy\in\scrP$ and $0\leq \ell <\tau_1(\vech,\hatp)$.
Hence, for $\hatq=(q_0,\vecq)$, $\hatp=(1,\vecp)$, the parameter $\ell$ describes the change in $q_0$ since the last collision. 

We now calculate the mean collision time {\em for fixed $\hatp$} (since $d\vecp$ is not a finite measure as opposed to the Lebesgue measure $d\vecv$ on $\S_1^{d-1}$) and proceed as in Section \ref{sec3}. The Lebesgue measure $d\hatq$ reads
\begin{equation}
d\hatq = d\vech\,d\ell .
\end{equation}
Repeating the steps that led to \eqref{mean1}, we find that the mean collision time is (denote by $\scrD_\hatp$ the projection of $\scrD$ onto $\{0\}\times\RR^n$ in direction $\hatp$)
\begin{equation}\label{mean100} 
\begin{split}
& \frac{1}{\#\scrP_T\vol(r\Sigma) + T^{d-1}\vol\scrD_\hatp} 
\times \\ 
& \times
\bigg\{ \sum_{\vecy\in\scrP_T} \int_{\{0\}\times r\Sigma}
\tau_1(\vech+\vecy, \hatp)\, d\vech +  \int_{T\scrD_\hatp} \tau_1(\vecq(\vech),\hatp) \, d\vech \bigg\}.
\end{split}
\end{equation}
This equals
\begin{equation}
\frac{\vol(T\scrD)}{\#\scrP_T\vol(r\Sigma)+ T^{d-1}\vol\scrD_\hatp} ,
\end{equation}
which, in view of \eqref{density}, converges to (as $T\to\infty$)
\begin{equation}
\frac{1}{\nbar \vol(r\Sigma)}
=\frac{1}{\nbar\,\sigmabar r^{d-1}} =\frac{\xibar}{r^{d-1}}.
\end{equation}
We see that the mean collision time is independent of the choice of $\hatp=(1,\vecp)$. The mean free path length for fixed $\vecp$ is thus 
\begin{equation}
\frac{\|\vecp\|\xibar}{r^{d-1}} ,
\end{equation}
which scales in the same way as the mean collision time provided $\|\vecp\|\asymp 1$ as $r\to 0$.

\section{The Boltzmann-Grad limit of kicked Hamiltonians}\label{sec8}

By the same reasoning as for the Lorentz gas (Section \ref{sec4}), the formulas for the mean free path length and collision time suggest to introduce the macroscopic coordinates 
\begin{equation}
\hatQ(t)= r^{d-1} \hatq(r^{-(d-1)} t) = (t,r^{d-1} \vecq(r^{-(d-1)} t) ) ,
\end{equation}
\begin{equation}
\hatP(t)= \hatp(r^{-(d-1)} t) = (1,r^{d-1} \vecp(r^{-(d-1)} t) ) .
\end{equation}
We can ask the same questions as in Section \ref{sec4} and seek a limit process as $r\to 0$ which is governed by a transport equation of the form
\begin{equation}\label{FPK2}
(\partial_t +\vecP\cdot\partial_\vecQ-\partial_\xi) f_t(\vecQ,\vecP,\xi,\vecomega) \\
= [\scrC f_t](\vecQ,\vecP,\xi,\vecomega)
\end{equation}
with a collision operator $\scrC$ as in Section \ref{sec4}. Note that $\|\vecP\|$ is now  {\em not} a constant of motion, as it is not preserved by the scattering map.

\section{Renormalisation of the transition kernel for the Lorentz gas}\label{sec9}

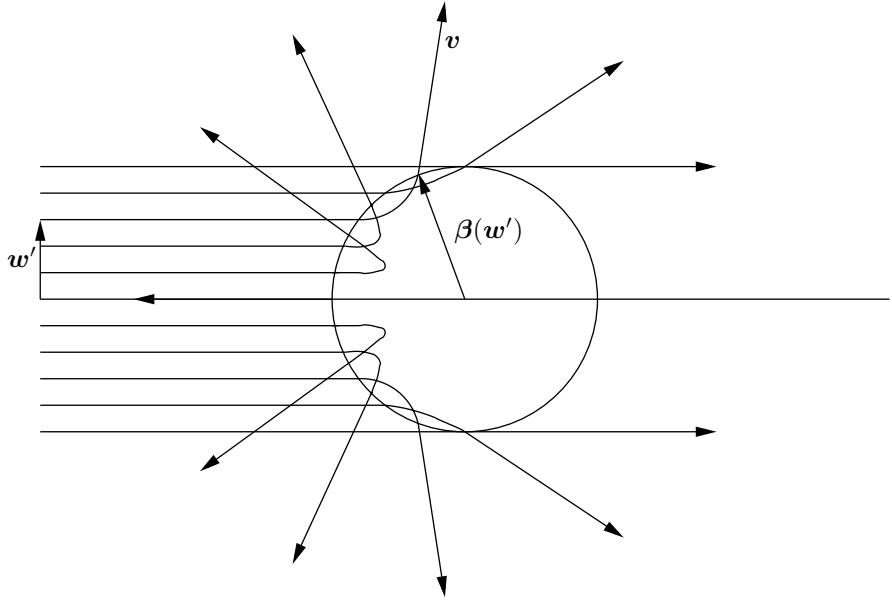
\begin{figure}
\begin{center}
\begin{minipage}{\textwidth}
\unitlength0.1\textwidth
\begin{picture}(10,6.3)(0,0)
\put(0.3,0){\includegraphics[width=0.9\textwidth]{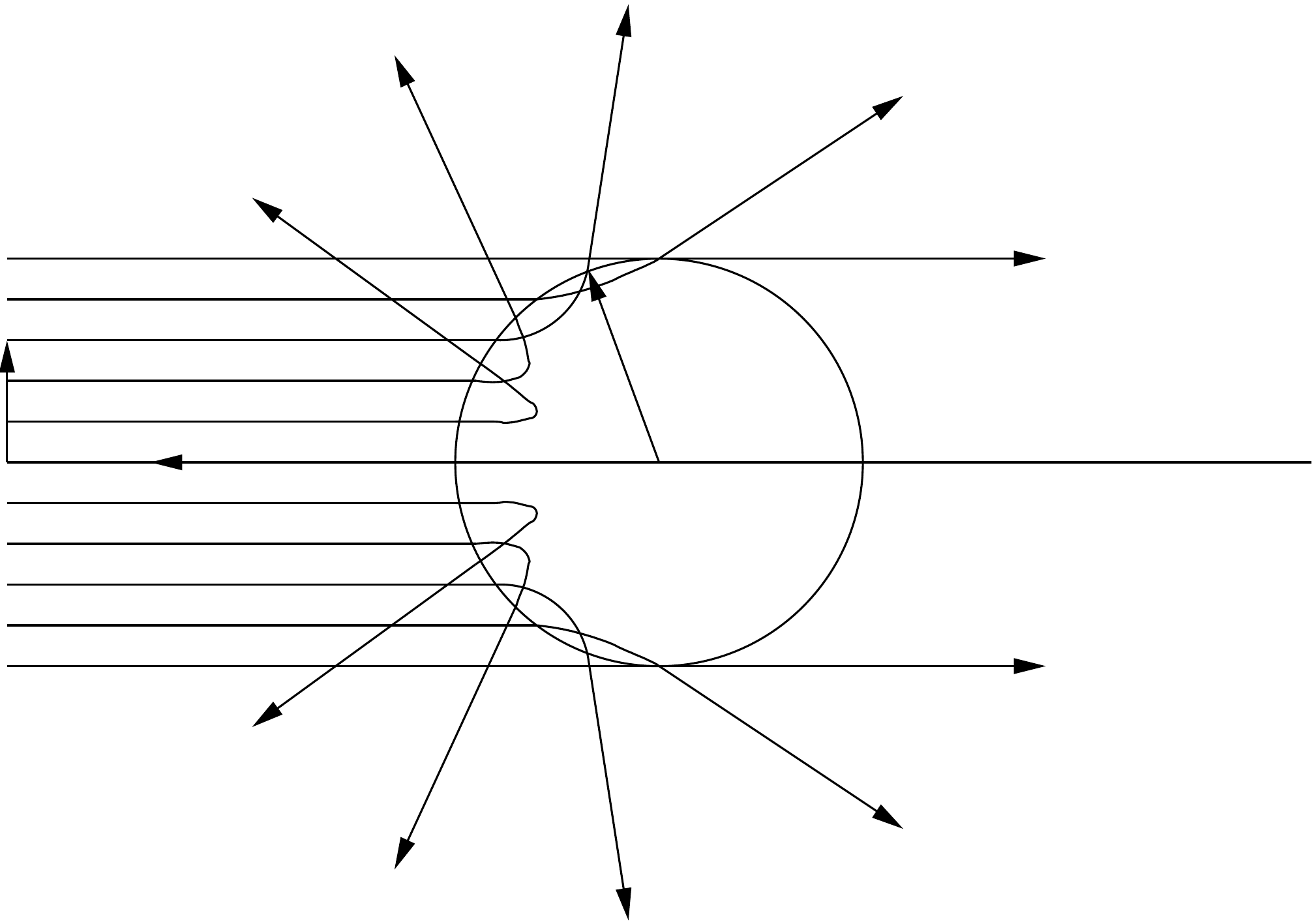}}
\put(0,3.5){$\vecw'$} 
\put(4.7,3.8){$\vecbeta(\vecw')$}
\put(4.6,5.8){$\vecv$} 
\end{picture}
\end{minipage}
\end{center}
\caption{Scattering of a parallel beam of particles.} \label{fig4}
\end{figure}

The plan is to understand how the transition kernel $k(\vecomega',\xi,\vecomega)$ emerges in the Boltzmann-Grad limit from the microscopic dynamics of the Lorentz gas. To this end consider a parallel beam of particles with velocity $\vecv'=\vece_1$ that hit a scatterer of the Lorentz gas with impact parameter $\vecw'\in\scrB_1^{d-1}$ as shown in Figure \ref{fig4}. 

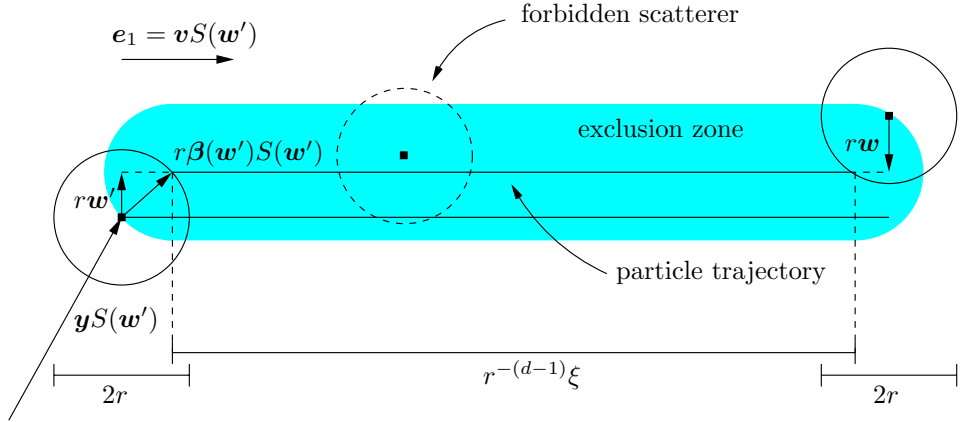
\begin{figure}
\begin{center}
\begin{minipage}{\textwidth}
\unitlength0.1\textwidth
\begin{picture}(10,4.5)(0,0)
\put(0,0){\includegraphics[width=\textwidth]{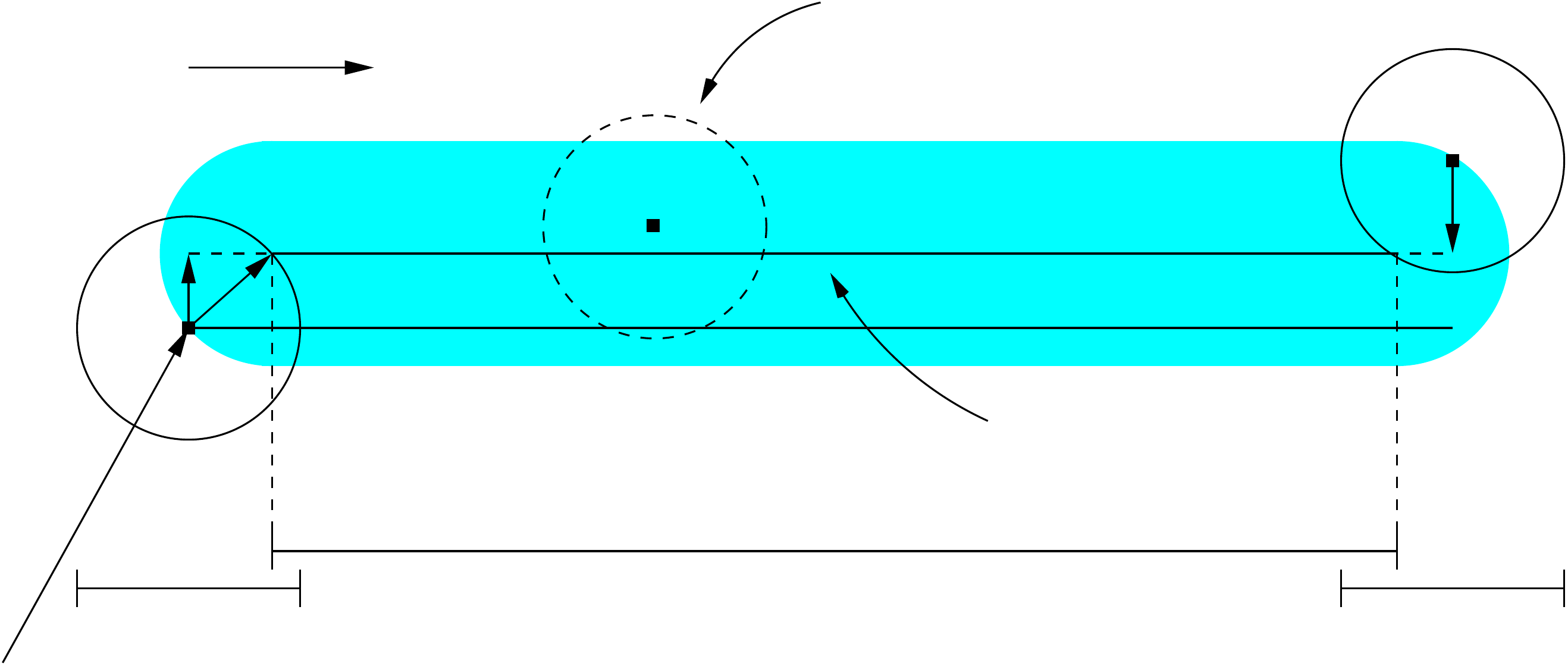}}
\put(0.7,2.25){$r\vecw'$} 
\put(1.75,2.75){$r\vecbeta(\vecw') S(\vecw')$} 
\put(0.7,1){$\vecy S(\vecw')$}
\put(1.1,4){$\vece_1=\vecv S(\vecw')$} 
\put(1,0.2){$2r$}
\put(9.1,0.2){$2r$}
\put(8.8,2.85){$r\vecw$}
\put(5,0.4){$r^{-(d-1)}\xi$}
\put(5.4,4.2){forbidden scatterer}
\put(6.4,1.5){particle trajectory}
\put(6,3){exclusion zone}
\end{picture}
\end{minipage}
\end{center}
\caption{Inter-collision flight in the Lorentz gas. The exclusion zone is the cylinder with caps described in the text.} \label{fig5}
\end{figure}

We assume $\vecw'$ is a random vector in $\scrB_1^{d-1}$, whose distribution is given by the probability measure $\lambda$ on $\scrB_1^{d-1}$ which is absolutely continuous with respect to Lebesgue measure. In general, $\lambda$ will depend on the entire history of the particle trajectories (and thus on $r$), which leads to significant technical complications that are resolved in \cite{Marklof:2011ho} in the case of the periodic Lorentz gas. For the purpose of this lecture, we simply assume that $\lambda$ is arbitary but fixed. The task is now to calculate the probability of hitting the next scatterer around time $r^{-(d-1)}\xi$ with impact parameter near $\vecw$. It is convenient to rotate our coordinate system by $S(\vecw')\in\SO(d)$, so that the outgoing velocity $\vecv$ becomes $\vece_1$. Recall from \eqref{inout2} that $\vecv =  \vece_1 S(\vecw')^{-1}$.

Let us denote by $\vecy$ the current scatterer location.
Define the open cylinder of length $\ell$ and radius $r$,
\begin{equation}
\fZ(\ell,r)=(0,\ell)\times\scrB_r^{d-1}
\end{equation}
and the corresponding cylinder with spherical caps,
\begin{equation}
\fZ_\text{caps}(\ell,r) = \fZ(\ell,r) \cup \scrB_r^d \cup (\ell \vece_1 +\scrB_r^{d}).
\end{equation}
The probability of hitting the next scatterer at a time in the interval $[\frac{\xi}{r^{d-1}},\frac{\xi+\epsilon}{r^{d-1}})$ with impact parameter in the ball $\vecw+\scrB_{r\delta}^{d-1}$ for some small $\epsilon>0$, $\delta>0$ is the probability (i.e. the $\lambda$-measure of the set of $\vecw'\in\scrB_1^{d-1}$) that 
\begin{enumerate}[(i)]
\item $$\bigg(\fZ_\text{caps}\bigg(\frac{\xi}{r^{d-1}},r\bigg) + r \big(\sqrt{1-{b'}^2}, \vecw'\big) \bigg) \cap (\scrP-\vecy) S(\vecw')= \emptyset,$$
\item 
\begin{multline*}
\bigg( \bigg(\frac{\xi}{r^{d-1}}+r\sqrt{1-{b'}^2}+r\sqrt{1-b^2}, r (\vecw'-\vecw) \bigg) \\
+ \bigg(\bigg[0,\frac{\epsilon}{r^{d-1}}\bigg) \times \scrB_{r\delta}^{d-1}\bigg) \bigg) \cap(\scrP-\vecy) S(\vecw') \neq \emptyset,
\end{multline*}
\end{enumerate}
cf.~Figure \ref{fig5}.
We assume here and in the following that the probability that there are two or more elements of $(\scrP-\vecy) S(\vecw')$ in a set of small diameter is small. (This can be proved for all examples discussed in these lectures.) 

Note that for
\begin{equation}
D(r) =
\begin{pmatrix}
r^{d-1} & \vecnull \\ \trans\vecnull  & r^{-1} 1_{d-1}
\end{pmatrix}
\end{equation}
we have
\begin{equation}
\fZ\bigg(\frac{\xi}{r^{d-1}},r\bigg)  D(r) = \fZ(\xi,1) 
\end{equation}
and
\begin{equation}
\fZ_\text{caps}\bigg(\frac{\xi}{r^{d-1}},r\bigg) D(r) \xrightarrow[r\to 0]{} \fZ(\xi,1) .
\end{equation}
Up to a small error (as $r\to 0$) we can therefore replace (i), (ii) by 
\begin{enumerate}[(i)]
\item $\fZ(\xi,1) \cap \Theta_r(\vecw') = \emptyset $,
\item $\big([0,\epsilon) \times \scrB_{\delta}^{d-1}\big) +(\xi,-\vecw) \cap \Theta_r(\vecw')\neq \emptyset$,
\end{enumerate}
where
\begin{equation}\label{theta1}
\Theta_r(\vecw') = (\scrP-\vecy) S(\vecw') D(r) - (0,\vecw')
\end{equation}
defines a sequence (in $r>0$) of random point processes for $\vecw'$ distributed according to $\lambda$. 

The main objective is to show that (for every fixed $\vecy\in\scrP$, or if this is not possible, for suitably random $\vecy$) there is a random process $\Theta$ such that
\begin{equation}\label{limitTheta}
\Theta_r(\vecw') \xrightarrow[r\to 0]{} \Theta - (0,\vecw')
\end{equation}
in distribution (we in fact only require convergence for two-dimensional distributions), where $\Theta$ and the random variable $\vecw'$ are independent. The limit process $\Theta$ can then be used to find a formula for the transition kernel $k(\vecomega',\xi,\vecomega)$ with 
\begin{equation}
\vecomega'=\vecomega'(\vecs',\ldots),\qquad \vecomega=\vecomega(\vecb,\ldots).
\end{equation}
The requirement to keep track of other ``hidden'' variables (``\ldots'') will depend on the particular choice of scatterer configuration $\scrP$, since the limit process $\Theta$ can in principle depend on $\vecy$. We will illustrate this strategy with a few natural examples in Sections \ref{sec11}--\ref{sec19}. Let us first discuss what the above argument looks like in the case of kicked Hamiltonians.

\begin{figure}
\begin{center}
\begin{minipage}{\textwidth}
\unitlength0.1\textwidth
\begin{picture}(10,4.5)(0,0)
\put(0,0){\includegraphics[width=\textwidth]{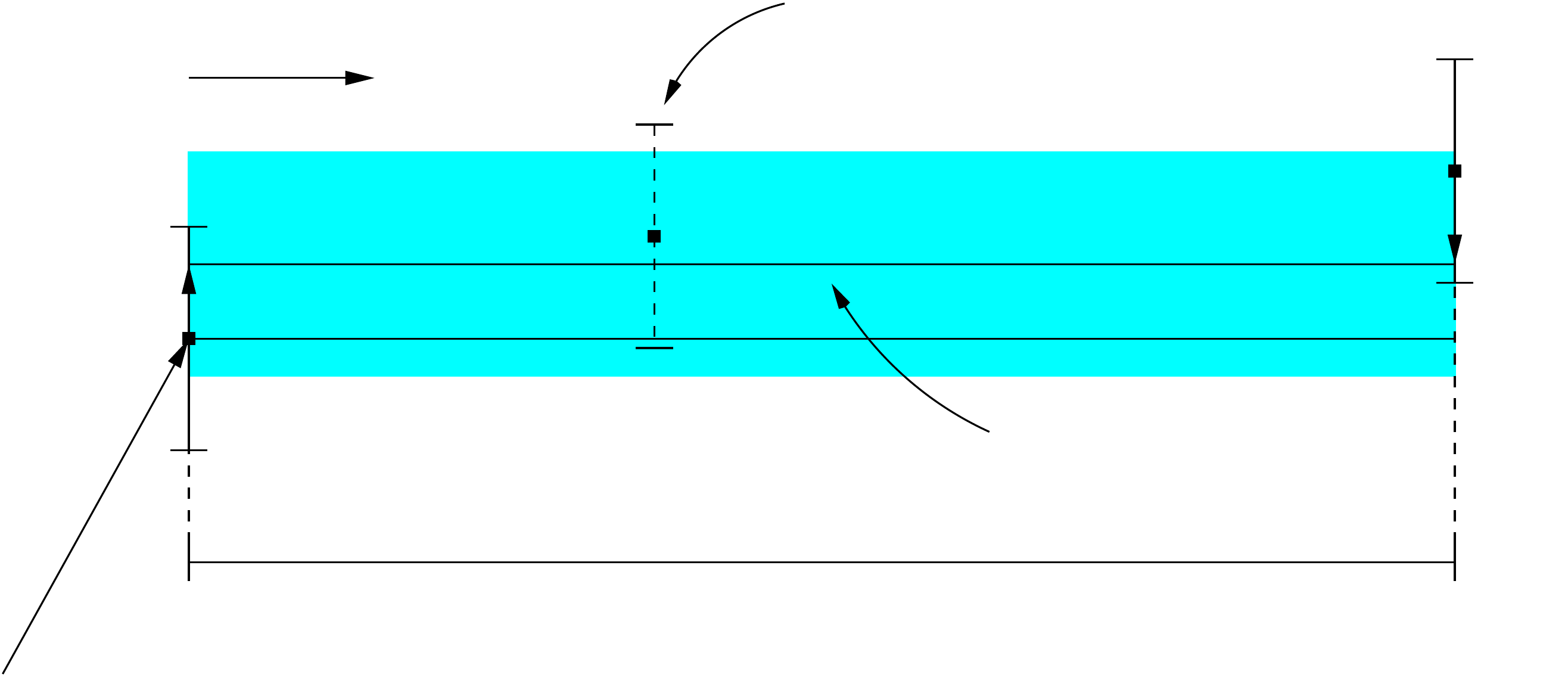}}
\put(0.7,2.25){$r\vecw'$} 
\put(0.3,0.3){$\vecy S(\vecw')$}
\put(1.1,4){$\vece_1=\hatp S(\vecw')$} 
\put(8.8,2.85){$r\vecw$}
\put(5,0.4){$r^{-(d-1)}\xi$}
\put(5.1,4.2){forbidden scatterer}
\put(6.4,1.5){particle trajectory}
\put(6,3){exclusion zone}
\end{picture}
\end{minipage}
\end{center}
\caption{Inter-collision flight for a kicked Hamiltonian. The exclusion zone is the cylinder without caps described in the text.} \label{fig6}
\end{figure}
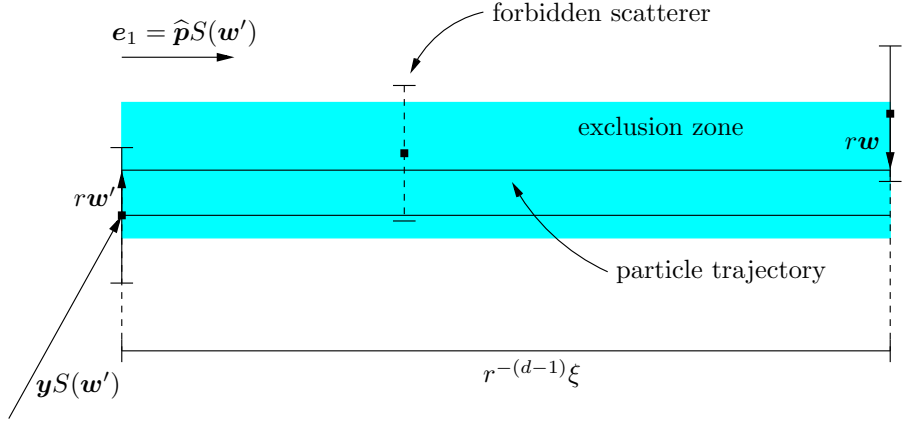

\section{Renormalisation of the transition kernel for kicked Hamiltonians}\label{sec10}

Recall that in this setting the scattering is described by the matrix $S(\vecw)$ as defined in \eqref{SbKH}, and the impact parameter $\vecw'$ is random according to an absolutely continuous probability measure $\lambda$ on the cross section $\Sigma$.
We now use the cylinder 
\begin{equation}
\fZ_\Sigma(\ell,r)=(0,\ell)\times r\Sigma
\end{equation}
to account for the different scattering cross section, see Figure \ref{fig6}. Following the same argument as in Section \ref{sec9}, we find that the probability of hitting the next scatterer at a time in $[\frac{\xi}{r^{d-1}},\frac{\xi+\epsilon}{r^{d-1}})$ with impact parameter in $\vecw+\scrB_{r\delta}^{d-1}$ equals (this time no approximation is needed)
\begin{enumerate}[(i)]
\item $\fZ_\Sigma(\xi,1) \cap \Theta_r(\vecw') = \emptyset $,
\item $\big([0,\epsilon) \times \scrB_{\delta}^{d-1}\big) +(\xi,-\vecw) \cap \Theta_r(\vecw')\neq \emptyset$,
\end{enumerate}
with the random process
\begin{equation}\label{theta2}
\Theta_r(\vecw') = (\scrP-\vecy) S(\vecw') D(r) - (0,\vecw') .
\end{equation}
The question is, as before, does 
\begin{equation}\label{limitKH}
\Theta_r(\vecw') \xrightarrow[r\to 0]{} \Theta - (0,\vecw')
\end{equation}
hold with $\Theta$ and $\vecw'$ independently distributed?

\section{Poisson process}\label{sec11}

The Poisson process $\Theta$ in $\RR^d$ with intensity $\nbar$ is characterised by the property that for any collection of bounded, pairwise disjoint Borel sets $\scrB_1,\ldots,\scrB_k$ and integers $r_1,\ldots,r_k\geq 0$ we have the probability
\begin{equation}
\PP(\# (\Theta\cap\scrB_1)=r_1,\ldots, \# (\Theta\cap\scrB_k)=r_k ) 
= \prod_{i=1}^k \frac{(\nbar\vol\scrB_i)^{r_i}}{r_i!} \e^{-\nbar\vol\scrB_i} .
\end{equation}
If the limit process $\Theta$ in Section \ref{sec9} or \ref{sec10} is a Poisson process (conditioned so that $\vecnull\in\Theta$), then the probability of the events (i) and (ii) is (set $\Sigma=\scrB_1^{d-1}$ in the case of the Lorentz gas)
\begin{multline}
\e^{-\nbar\vol\fZ_\Sigma(\xi,1)} (1- \e^{-\nbar\epsilon \vol\scrB_\delta^{d-1}}) \\
=
\e^{-\nbar\vol\fZ_\Sigma(\xi,1)} \nbar\epsilon \vol\scrB_\delta^{d-1} + \text{(lower order terms in $\epsilon,\delta$)}
\end{multline}
and thus, with $\epsilon\to d\xi$, $\vol\scrB_\delta^{d-1}\to d\vecb$, we find for the limit transition kernel
\begin{equation}
k(\vecomega',\xi,\vecomega)=\xibar^{-1} \; \e^{-\xi/\xibar} , \qquad \xibar:=\frac{1}{\nbar\,\sigmabar},  
\end{equation}
with 
\begin{equation}
\vecomega:= \vecb , \qquad d\pp(\vecomega) :=  \sigmabar^{-1} d\vecb
\end{equation}
and the total scattering cross section $\sigmabar=\vol\Sigma$.
This is the same kernel as in the example \eqref{ex0} leading to the linear Boltzmann equation.

It follows from a standard probabilistic argument (cf.~\cite{Boldrighini:1983jm}) that, if $\scrP$ is a typical realisation of a Poisson process, then $\Theta_r\to\Theta-(0,\vecw')$ with $\Theta$ Poisson (again conditioned to $\vecnull\in\Theta$). Let us now discuss examples of scatterer configurations, for which $\Theta$ is {\em not} Poisson and for which the linear Boltzmann equation fails.

\section{Point processes and homogeneous spaces}\label{sec12}

In the following we will consider homogeneous spaces of the form 
\begin{equation}
\GamG=\{ \Gamma g : g\in G\}
\end{equation}
where $G$ is a connected Lie group and $\Gamma$ is a {\em lattice} in $G$. {\em Lattice} means $\Gamma$ is a discrete subgroup of $G$ such that there exists a fundamental domain $\scrF_\Gamma$ of $\Gamma$ in $G$ with finite left Haar measure $\mu$. 

It is a general property of the left Haar measure that it is, up to a character, also right invariant. That is, there exists a character $\chi: G\to \RR$ such that, for any measurable set $\scrB$ and any $g\in G$, we have $\mu(\scrB g)=\chi(g)\mu(\scrB)$. It is an easy exercise to show that every fundamental domain of $\Gamma$ in $G$ has the same finite left Haar measure. Since $\scrF_\Gamma g$ is a fundamental domain for every $g\in G$, we see that $\mu(\scrF_\Gamma g)=\mu(\scrF_\Gamma)$ and hence $\chi=1$. This means $G$ is unimodular, i.e., left and right Haar measure coincide. 

We normalise the Haar measure so that $\mu(\scrF_\Gamma)=1$, and also denote by $\mu$ the push-forward of Haar measure to $\GamG$, which is the unique $G$-invariant (under right multiplication) probability measure on $\GamG$.

The plan is now to consider point patterns $\scrP$, whose corresponding point processes $\Theta_r$ and $\Theta$ can be parametrised by a homogeneous space of the above type. The convergence of $\Theta_r$ will follow from the ergodic theory for subgroups of $G$ acting on $\GamG$ by right translation.

We will start developing our theory by starting with the simplest example for $\GamG$: the space of lattices.

\section{The space of lattices}\label{sec13}

The opposite extreme of a random scatterer configuration is a perfectly periodic point set $\scrP$. We begin by assuming that $\scrP$ is a Euclidean lattice $\scrL$. We fix the covolume of $\scrL$ (the volume of its fundamental domain) to be $\nbar^{-1}$, so that the asymptotic density of $\scrP=\scrL$ is $\nbar$, as assumed in Section \ref{sec2}. Every such lattice can be written as
\begin{equation}
\scrL = \nbar^{-1/d}\ZZ^d M
\end{equation}
for some $M\in\SLdR$. Since the stabiliser of $\ZZ^d$ under right multiplication by $G=\SLdR$ is the subgroup $\Gamma=\SLdZ$, one can show that there is a bijection
\begin{equation}\label{SOL}
\begin{split}
\GamG & \xrightarrow{\sim} \{ \text{Euclidean lattices of covolume $\nbar^{-1}$} \} \\
\Gamma M & \mapsto \nbar^{-1/d}\ZZ^d M .
\end{split}
\end{equation}
To find the inverse map, note that, for any basis $\veca_1,\ldots,\veca_d$ of $\scrL$, the matrix $M=\nbar^{1/d} (\trans\veca_1,\ldots,\trans\veca_d)$ is in $\SLdR$; the substitution $M\mapsto \gamma M$, $\gamma\in\Gamma$ corresponds to a base change of $\scrL$. 

It is a well known fact (due to Minkowski) that $\Gamma=\SLdZ$ is a lattice in $G=\SLdR$. We can use the $G$-invariant probability measure $\mu$ on $\GamG$ to define a random point process in $\RR^d$ by setting 
\begin{equation}\label{ThetaSOL}
\Theta = \nbar^{-1/d}\ZZ^d M
\end{equation}
with $M$ random in $\GamG$ according to $\mu$ and the above identification \eqref{SOL} of $\GamG$ and the space of lattices. Alternatively, we could also take $M$ in \eqref{ThetaSOL} to be random according to Haar measure in a fundamental domain $\scrF_\Gamma$ in $G$.

The following theorem says that, for any fixed $\scrP=\scrL$ as above,
\begin{equation}
\Theta_r \to \Theta - (0,\vecw')
\end{equation}
where $\Theta$ and $\vecw'$ are independent, as required.

\begin{thm}[\cite{Marklof:2010ib}]\label{thm13.1}
Let $\lambda$ be an absolutely continuous probability measure on $\Sigma$, $\scrA_1,\ldots,\scrA_k\subset\RR^d$ bounded with boundary of measure zero and $r_1,\ldots,r_k\in\ZZ_{\geq 0}$. Then
\begin{equation}\label{statem}
\lambda(\#(\Theta_r\cap\scrA_i) = r_i)  \xrightarrow[r\to 0]{} 
\int_{\Sigma} \mu\big( \#((\Theta-(0,\vecw'))\cap \scrA_i)=r_i   \big) \, d\lambda(\vecw').
\end{equation}
\end{thm}

This theorem is in turn a consequence of the following equidistribution theorem on $\GamG$. The matrix $S(\vecw)$ is as defined in \eqref{SbLG} or \eqref{SbKH}, respectively. (In the former, take $\Sigma=\scrB_1^{d-1}$.)

\begin{thm}[\cite{Marklof:2010ib}]\label{thm13.2}
For any $M\in\GamG$, any bounded continuous $f: \Sigma\times \GamG \to \RR$ and any absolutely continuous probability measure $\lambda$ on $\Sigma$, 
 \begin{equation}\label{eq13.2i}
\int_{\Sigma} f(\vecw', M S(\vecw') D(r))\, d\lambda(\vecw')  \xrightarrow[r\to 0]{}
\int_{\Sigma} \int_{\GamG} f(\vecw',M) \,d\mu(M)\, d\lambda(\vecw').
\end{equation}
\end{thm}
This states the equidistribution of large spheres (for $S(\vecw)$ as in \eqref{SbLG}), or expanding horospheres (for $S(\vecw)$ as in \eqref{SbKH}). 

Statement \eqref{statem} in Theorem \ref{thm13.1} follows from Theorem \ref{thm13.2} by choosing in \eqref{eq13.2i} as test function $f$ the characteristic function of the set
\begin{equation}\label{0set}
\big\{ (\vecw',M) \in \Sigma \times \GamG :  \#\big((\nbar^{-1/d}\ZZ^d M - (0,\vecw' )\cap\scrA_i \big)= r_i \big\} .
\end{equation}
Now, this characteristic function is of course not continuous, but one can show that \eqref{0set} has boundary of measure zero, and thus the characteristic function can be approximated sufficiently well by continuous functions. Details of this technical argument can be found in \cite{Marklof:2010ib}, Sections 5 and 6.

We will later (in Section \ref{sec18}) return to the justification of statements of the above type, in fact, in a more general form (Theorem \ref{thm18.1}). Let us first see how Theorem \ref{thm13.1} yields an expression for the transition kernel $k(\vecomega',\xi,\vecomega)$.

\section{The transition kernel for lattices}\label{sec14}

To complete the programme of Section \ref{sec9} (for the Lorentz gas) and Section \ref{sec10} (for kicked Hamiltonians) we apply Theorem \ref{thm13.1} for $k=2$, $r_1=0$, $r_2=1$ and the test sets
\begin{equation}
\scrA_1=\fZ_\Sigma(\xi,1),\qquad \scrA_2= ([0,\epsilon)\times\scrB_\delta^{d-1}) +(\xi,-\vecw)
\end{equation}
with $\epsilon,\delta$ small.

Note that the condition 
\begin{equation}
\#\big( (\Theta-(0,\vecw'))\cap\scrA_2 \big) =1
\end{equation}
means that
\begin{equation}
\#\big( (\Theta-(\xi,\vecw'-\vecw))\cap([0,\epsilon)\times\scrB_\delta^{d-1}) \big) =1,
\end{equation}
i.e., our random lattice $\Theta$ must have a lattice point near $(\xi,\vecw'-\vecw))$. In order to characterise lattices with this property, set $X=\GamG$ and define the subspace
\begin{equation}
X(\vecy)=\{ M\in X : \vecy\in\ZZ^d M \}
\end{equation}
for a given $\vecy\in\RR^d$. In other words, $X(\vecy)$ is the subspace of all lattices (of covolume one) that contain $\vecy$.
In \cite{Marklof:2010ib} we construct a probability measure $\nu_\vecy$ on $X(\vecy)$ so that 
\begin{equation}\label{decoy}
d\mu(M) = d\nu_\vecy(M)\, d\vecy.
\end{equation}
(Since the spaces $X(\vecy)$ are not disjoint, formula \eqref{decoy} is only valid for Borel sets $\scrE\subset  X$, $U\subset\RR^d\setminus\{\vecnull\}$ so that for all $\vecy_1\neq\vecy_2\in U$ we have $X(\vecy_1)\cap X(\vecy_2)\cap\scrE=\emptyset$. See \cite[Prop.~7.3]{Marklof:2010ib} for details.) 

The decomposition \eqref{decoy} shows that the transition kernel is, with $\vecomega'=\vecs'=(0,\vecw')$, $\vecomega=\vecb=(0,\vecw)$, $d\pp(\vecomega)=\sigmabar^{-1} d\vecw$ and $\xibar=\frac{1}{\nbar\,\sigmabar}$, 
\begin{equation}\label{1st}
k(\vecomega',\xi,\vecomega) = \xibar^{-1} \nu_\vecy\big(\big\{  M\in X(\vecy) : \nbar^{-1/d} \ZZ^d M \cap (\fZ_\Sigma(\xi,1)+(0,\vecw'))=\emptyset \big\}\big)
\end{equation}
where $\vecy=\nbar^{1/d} (\xi,\vecw'-\vecw)$. For an explicit description of the $\nu_\vecy$-measure of the above set, see \cite{Marklof:2011di}, Section 2.2.

We note that $\nu_\vecy$ is invariant under $G$, in the sense that for any Borel set $\scrE\subset X(\vecy)$ we have for all $g\in G$ 
\begin{equation}
\nu_\vecy(\scrE)=\nu_{\vecy g}(\scrE g).
\end{equation}
In particular for
\begin{equation}
g = \begin{pmatrix} \nbar^{1-1/d} & \vecnull \\ \trans\vecnull & \nbar^{-1/d} 1_{d-1} \end{pmatrix} \in G
\end{equation}
we see that $\vecy g= (\nbar\xi,\vecw'-\vecw)$,
and furthermore
\begin{equation}
\nbar^{-1/d} \ZZ^d M \cap (\fZ_\Sigma(\xi,1)+(0,\vecw')) g =\emptyset
\end{equation}
is equivalent to
\begin{equation}
\ZZ^d M \cap (\fZ_\Sigma(\nbar\xi,1)+(0,\vecw')) =\emptyset .
\end{equation}
Hence we obtain the formula 
\begin{equation}\label{1st1}
k(\vecomega',\xi,\vecomega) = \xibar^{-1} \nu_\vecy\big(\big\{  M\in X(\vecy) : \ZZ^d M \cap (\fZ_\Sigma(\nbar \xi,1)+(0,\vecw'))=\emptyset \big\}\big)
\end{equation}
with $\vecy=(\nbar\xi,\vecw'-\vecw)$. The above scaling property must of course hold a priori, as it corresponds to a simple rescaling of length units.

In dimension $d=2$, when $\Sigma=\scrB_1^1=(-1,1)$, eq.~\eqref{1st1} can be used to calculate an explicit formula for the transition kernel. We have \cite{Marklof:2008dr} (with $\sigmabar=2$)
\begin{equation}\label{twodidi}
k(\vecomega',\xi,\vecomega)=  \frac{12\,\nbar}{\pi^2}\, \Upsilon\bigg( 1+ \frac{(\nbar\xi)^{-1} - \max(|\vecw|,|\vecw'|)-1}{|\vecw-\vecw'|}\bigg)
\end{equation}
with
\begin{equation}
\Upsilon(x)=
\begin{cases} 
0 & \text{if }x\leq 0\\
x & \text{if }0<x<1\\
1 & \text{if }1\leq x.
\end{cases}
\end{equation}
For independent derivations of Formula \eqref{twodidi} that do not employ eq.~\eqref{1st1}, see the papers by Caglioti and Golse \cite{Caglioti:2008um,Caglioti:2010ir} and Bykovskii and Ustinov \cite{Bykovskiui:2009db}.

There are no such formulas in higher dimension, although \eqref{1st1} can be used to extract information to obtain asymptotics for $\xi\to 0$ and $\xi\to\infty$, cf.~\cite{Marklof:2011di}. We have in particular in the case $\Sigma=\scrB_1^{d-1}$ 
\begin{equation}\label{Kbound}
\frac{1-2^{d-1} \xibar^{-1}\xi}{\zeta(d)\xibar} \leq k(\vecomega',\xi,\vecomega) \leq \frac{1}{\zeta(d)\xibar},
\end{equation}
and so for small $\xi$ this implies $k(\vecomega',\xi,\vecomega) = \zeta(d)^{-1} \xibar^{-1}+O(\xi)$.
Here $\zeta(d)$ is the Riemann zeta function and $\zeta(d)^{-1}$ is the relative density of primitive lattice points in $\ZZ^d$. (The {\em primitive} lattice points are those points in $\ZZ^d$ that are visible from the origin.) 
Compare \eqref{Kbound} with the result for the Poisson process (Section \ref{sec11}):
\begin{equation}
k_\text{Poisson}(\vecomega',\xi,\vecomega) = \xibar^{-1} \e^{-\xi/\xibar} = \xibar^{-1} - \xibar^{-2} \xi + O(\xi^2).
\end{equation}
The asymptotics of $k(\vecomega',\xi,\vecomega)$ for large $\xi$ is more complicated to state (see \cite{Marklof:2011di} for what we know). In the next section we will provide tail estimates for the distribution of free path lengths. 

To check the symmetry property $k(\vecomega,\xi,\vecomega')=k(\vecomega',\xi,\vecomega)$ in \eqref{ksymm}, note that in \eqref{1st1} we can replace the term $(0,\vecw')$ by $(0,\vecw')-\vecy=(-\nbar\xi,\vecw)$ since $\vecy\in  \ZZ^d M$. This, upon reflecting all sets at the origin, yields 
\begin{equation}
k(\vecomega',\xi,\vecomega) = \xibar^{-1} \nu_\vecy\big(\big\{  M\in X(\vecy) : \\ \ZZ^d M \cap (\fZ_{-\Sigma}(\nbar\xi,1)+(0,-\vecw))=\emptyset \big\}\big).
\end{equation}
Reflecting at the ``horizontal'' axis $\RR\vece_1$ yields
\begin{equation}
\vecy=(\xi,\vecw'-\vecw) \mapsto (\xi,\vecw-\vecw') ,
\end{equation}
\begin{equation}
\fZ_{-\Sigma}(\nbar\xi,1)\mapsto \fZ_\Sigma(\nbar\xi,1)
\end{equation}
and
\begin{equation}
(0,-\vecw)\mapsto (0,\vecw).
\end{equation}
The measure $\nu_\vecy$ is preserved under this map and hence we see that indeed $k(\vecomega,\xi,\vecomega')=k(\vecomega',\xi,\vecomega)$.

\section{The distribution of free path lengths}\label{sec15}

The limit distribution of the free path lengths (between consecutive collisions) can be obtained from the transition kernel via the formula \eqref{FreePL}.
As for the transition kernel, explicit formulas are only known in two dimensions, and were first computed heuristically by Dahlqvist \cite{Dahlqvist:1997bd}, with subsequent rigorous proof by Boca and Zaharescu \cite{Boca:2007kw}. These formulas can also be derived via \eqref{FreePL} from \eqref{twodidi}, see \cite{Marklof:2008dr}.
For arbitrary dimension $d\geq 2$, the expression \eqref{1st} is used in \cite{Marklof:2011di} to obtain the following tail asymptotics for $\Sigma=\scrB_1^{d-1}$:
\begin{equation}
\overline\Phi_0(\xi) = \frac{1}{\xibar\zeta(d)} +O(\xi) \qquad (\xi>0),
\end{equation}
\begin{equation}\label{large}
\overline\Phi_0(\xi) \sim \frac{C}{\xi^3} \qquad (\xi\to\infty)
\end{equation}
with the constant
\begin{equation}
C=\frac{2^{2-d}}{d(d+1)\nbar^2\zeta(d)}.
\end{equation}
These asymptotics sharpen earlier upper and lower bounds by Bourgain, Golse and Wennberg \cite{Bourgain:1998gu,Golse:2000ka}.
Note that \eqref{large} implies that the density $\overline\Phi_0(\xi)$ has no second moment. This fact is used in \cite{super} to prove a superdiffusive central limit theorem for the periodic Lorentz gas. The first moment is of course, by our normalisation, the mean free path length
\begin{equation}
\int_0^\infty \xi \overline\Phi_0(\xi)\,d\xi =\xibar.
\end{equation}
Compare these results with the Poisson case:
\begin{equation}
\overline\Phi_{0,\text{Poisson}}(\xi) = \xibar^{-1} \e^{-\xi/\xibar} .
\end{equation}
The main difference of the two is clearly the exponential vs.\ power-law tail.

\section{Cut-and-project sets}\label{sec16}

We will now significantly generalise the above discussion by allowing $\scrP$ to be a {\em cut-and-project set}, following closely the presentation in \cite{Marklof2013}. Examples include the classic quasicrystals (such as the vertex set of a Penrose tiling), and locally finite periodic point sets. A cut-and-project set $\scrP\subset\RR^d$ is defined as follows. 

For $m\geq 0$, $n=d+m$, let
\begin{equation}
\pi : \RR^{n} \to \RR^d ,
\end{equation}
\begin{equation}
\pi_\intl : \RR^{n} \to \RR^m 
\end{equation}
be the orthogonal projections of $\RR^{n}=\RR^d\times\RR^m$ onto the first and second factor, respectively. $\RR^d$ will be called the {\em physical space}, and $\RR^m$ the {\em internal space}. Let $\scrL\subset\RR^{n}$ be a lattice of full rank. Then the closure
\begin{equation}
\scrA := \overline{\pi_\intl(\scrL)} \subset\RR^m
\end{equation}
is an abelian subgroup. We denote by $\scrA^0$ the connected component of $\scrA$ containing $0$. $\scrA^0$ is a linear subspace of $\RR^m$ of dimension $m_1$. We find vectors $\veca_1,\ldots,\veca_{m_2}$ ($m=m_1+m_2$) so that
\begin{equation}
\scrA = \scrA^0 \oplus \ZZ\pi(\veca_1)\oplus \ldots \oplus \ZZ\pi(\veca_{m_2}) .
\end{equation}
The Haar measure of $\scrA$ is denoted by $\mu_\scrA$ and normalised so that $\mu_\scrA\big|_{\scrA^0}$ is the standard Lebesgue measure on $\scrA^0$.

For $\scrV:=\RR^d\times \scrA^0$, we note that $\scrL\cap\scrV$ is a full rank lattice in $\scrV$. For $\scrW\subset\scrA$ with non-empty interior, we call 
\begin{equation}
\scrP=\scrP(\scrW,\scrL)=\{ \pi(\vecell) : \vecell\in\scrL, \; \pi_\intl(\vecell)\in\scrW \}
\end{equation}
a {\em cut-and-project set}. $\scrW$ is called the {\em window set}. If the boundary of the window set has $\mu_\scrA$-measure zero, we say $\scrP(\scrW,\scrL)$ is {\em regular}. We will furthermore assume that $\scrW$ and $\scrL$ are chosen so that the map
\begin{equation}
\pi_\scrW : \{ \vecell\in\scrL : \pi_\intl(\vecell)\in\scrW\} \to \scrP
\end{equation}
is bijective. This is to avoid coincidences in $\scrP$.

It follows from Weyl equidistribution that such $\scrP$ have density
\begin{equation}
\nbar = \frac{\mu_\scrA(\scrW)}{\vol(\scrV/(\scrL\cap\scrV))} .
\end{equation}

In the case of lattices, we were able to assume, by the translational symmetry of the lattice, that the scatterer location $\vecy$ from which we launch our particle is without loss of generality $\vecy=\vecnull$. This is of course not the case for a general cut-and-project set considered now. Note however that for $\vecy\in\scrP$ there is $\vecell\in\scrL$ such that $\vecell=\pi(\vecy)$ and
\begin{equation}\label{eq16.8}
\scrP(\scrW,\scrL) -\vecy = \scrP(\scrW-\vecy_\intl,\scrL) , \qquad \vecy_\intl:=\pi_\intl(\vecell).
\end{equation}
By adjusting $\scrW$ in this way, we may therefore assume in the following that $\vecy=\vecnull$; but keep in mind that the limit process depends on the choice of $\scrW$ and therefore on the scatterer location $\vecy$.

\section{Spaces of cut-and-project sets}\label{sec17}

The aim is now to describe the ``closure'' (in a suitable sense) of the orbit of $\scrP$ under the $\SLdR$-action and construct a probability measure on it. This will yield, as we shall see, our limit random process $\Theta$.

Set $G=\SL(n,\RR)$, $\Gamma=\SL(n,\ZZ)$ and define the embedding (for any $g\in G$) 
\begin{equation}
\begin{split}
\varphi_g  : \SLdR & \hookrightarrow G \\
A & \mapsto g \begin{pmatrix} A & 0_{d\times m} \\ 0_{m\times d} & 1_m \end{pmatrix} g^{-1} .
\end{split}
\end{equation}
Since $\SLdR$ is generated by unipotent subgroups, Ratner's theorems \cite{Ratner:1991dc} imply that there is a (unique) closed connected subgroup $H_g\leq G$ such that:
\begin{enumerate}[(i)]
\item $\Gamma\cap H_g$ is a lattice in $H_g$.
\item $\varphi_g(\SLdR)\subset H_g$.
\item The closure of $\Gamma\backslash\Gamma\varphi_g(\SLdR)$ is $\Gamma\backslash\Gamma H_g$.
\end{enumerate}
We will call $H_g$ a {\em Ratner subgroup.} We denote the unique $H_g$-invariant probability measure on $\Gamma\backslash\Gamma H_g$ by $\mu_{H_g}=\mu_g$.
Note that $\Gamma\backslash\Gamma H_g$ is isomorphic to the homogeneous space $(\Gamma\cap H_g)\backslash H_g$.

Pick $g\in G$, $\delta>0$ such that $\scrL=\delta^{1/n} \ZZ^n g$. Then one can show \cite[Prop.~3.5]{Marklof2013} that 
\begin{equation}
\pi_\intl(\delta^{1/n} \ZZ^n h g) \subset \scrA \quad \text{for all $h\in H_g$,}
\end{equation}
\begin{equation}
\pi_\intl(\delta^{1/n} \ZZ^n h g) = \scrA \quad \text{for $\mu_g$-almost all $h\in H_g$.}
\end{equation}

The image of the map
\begin{equation}
\begin{split}
\Gamma\backslash \Gamma H_g & \to \{\text{point sets in $\RR^d$}\} \\
h & \mapsto \scrP(\scrW, \delta^{1/n} \ZZ^n h g)
\end{split}
\end{equation}
defines a space of cut-and-project sets, and the push-forward of $\mu_g$ equips it with a probability measure. We have thus defined a random point process $\Theta$ in $\RR^d$, which is $\SLdR$ invariant, and whose typical realisation is a cut-and-project set with window $\scrW$ and internal space $\scrA$. As we will see, this process is precisely the limit process $\Theta$ we are looking for.

\section{Equidistribution}\label{sec18}

The following equidistribution theorems generalise Theorem \ref{thm13.2} stated earlier. They are a consequence of Ratner's measure classification theorems \cite{Ratner:1991dc}, and in particular follow from a theorem of Shah \cite[Thm.~1.4]{Shah1996} on the equidistribution of translates of unipotent orbits. 

Recall that the matrix $S(\vecw)$ is as defined in \eqref{SbLG} or \eqref{SbKH}, respectively.

\begin{thm}[\cite{Marklof2013}]\label{thm18.1}
Fix $g\in G$.
For any bounded continuous $f: \Sigma \times \Gamma\backslash \Gamma H_g \to \RR$ and any absolutely continuous probability measure $\lambda$ on $\Sigma$, 
 \begin{equation}\label{eq18.2i}
\int_{\Sigma} f(\vecw', \varphi_g(S(\vecw') D(r)))\, d\lambda(\vecw') \xrightarrow[r\to 0]{}
\int_{\Sigma} \int_{\Gamma\backslash \Gamma H_g} f(\vecw',h) \,d\mu_g(h)\, d\lambda(\vecw').
\end{equation}
\end{thm}

Note that this theorem reduces to the statement of Theorem \ref{thm13.2} when $n=d$.
As in the case of the space of lattices (Section \ref{sec13}), Theorem \ref{thm18.1} is the key tool to prove:

\begin{thm}[\cite{Marklof2013}]\label{thm18.2}
If $\Theta_r$ is the sequence of random point processes (as defined in Section \ref{sec9} resp.~\ref{sec10}) corresponding to a regular cut-and-project set $\scrP=\scrP(\scrL,\scrW)$ and scatterer location $\vecy\in\scrP$, then
\begin{equation}\label{conny}
\Theta_r \xrightarrow[r\to 0]{} \Theta-(0,\vecw')
\end{equation}
where $\Theta$ is the random point process for the cut-and-project set $\scrP(\scrL,\scrW-\vecy_\intl)$  (as constructed in Section \ref{sec17}) and $\vecw'$ randomly distributed according to $\lambda$. 
\end{thm}

The convergence in \eqref{conny} is understood in the sense of Theorem \ref{thm13.1}, i.e., all test sets $\scrA_i$ are assumed to be bounded and have boundary of measure zero. Theorem \ref{thm18.2} implies that the limit transition kernel $k(\vecomega',\xi,\vecomega)$ is a function of 
\begin{equation}
\vecomega'=(\vecs',\vecy_\intl), \qquad \vecomega=(\vecb,\veceta_\intl),
\end{equation}
where $\vecy_\intl$ and $\veceta_\intl$ are the ``internal'' coordinates of the scatterers $\veceta,\vecy$ involved in the consecutive collisions, as defined \eqref{eq16.8}. For further details see the forthcoming paper \cite{Marklof2014}.

Let us now sketch the main ideas in the proof of Theorem \ref{thm18.2}. By a standard approximation argument (see \cite[p.~1973]{Marklof:2010ib} for details), where we localise $\lambda$ near any given fixed $\vecw_0$, we may reduce the proof to test functions $f$ that are independent on the first variable $\vecw'$. We thus need to establish
\begin{equation}\label{eq18.2ii}
\int_{\Sigma} f(\varphi_g(S(\vecw') D(r)))\, d\lambda(\vecw') \xrightarrow[r\to 0]{}
\int_{\Gamma\backslash \Gamma H_g} f(h) \,d\mu_g(h) .
\end{equation}
Let us take $S(\vecw)$ as defined in \eqref{SbKH}, i.e., 
\begin{equation}\label{SbKH2}
S(\vecw)=n(\partial W(\vecw)) , \qquad n(\vecx)=\begin{pmatrix} 1 & \vecx \\ \trans\vecnull & 1_{d-1} \end{pmatrix}.
\end{equation}
The case of $S(\vecw)$ as defined in \eqref{SbLG} is analogous and discussed in detail in \cite{Marklof:2010ib,Marklof2013}.
The variable substitution (recall \eqref{KHmap})
\begin{equation}\label{KHmap2}
\vecw' \mapsto \vecx=-\partial W(\vecw')
\end{equation}
yields for the left hand side of \eqref{eq18.2ii},
\begin{equation}\label{lolo}
\int_{\Sigma} f(\varphi_g(S(\vecw') D(r)))\, d\lambda(\vecw')
= \int f(\varphi_g(n(\vecx) D(r)))\, d\tilde\lambda(\vecx) ,
\end{equation}
where $\tilde\lambda$ is still absolutely continuous in view of \eqref{sigKH} and our assumption on the invertibility of the scattering map. Eq.~\eqref{lolo} defines a sequence of Borel probability measures $\rho_r$ on $\GamG$ via the linear functional
\begin{equation}
\rho_r[f] := \int f(\varphi_g(n(\vecx) D(r)))\, d\tilde\lambda(\vecx) .
\end{equation}
This sequence can be shown to be tight, which implies that any sequence of $\rho_r$ contains a convergent subsequence $\rho_{r_i}$. The question is: What are the possible weak limits $\rho_{r_i}\to\rho$? The first crucial observation is that, because $\varphi_g(n(\vecx) D(r)$ is a sequence of {\em expanding} horospheres, any limit must be invariant under the horospherical subgroup. That is,
\begin{equation}
\rho[f \circ \varphi_g(n(\vecx)) ] = \rho[f]
\end{equation}
for all $\vecx\in\RR^d$. Ratner's theorem \cite{Ratner:1991dc} gives a complete classification of such measures (in a significantly more general setting). In particular, the ergodic components are unique $H$-invariant probability measures $\mu_H$ supported on embedded homogeneous spaces of the form $\Gamma\backslash\Gamma H$ where $H$ is a closed connected subgroup of $G$ such that $\Gamma\cap H$ is a lattice in $H$. Shah \cite{Shah1996} showed (again in a more general setting) that for expanding translates any limit $\rho$ is in fact invariant under the group $\varphi_g(\SLdR)$, and $\rho=\mu_H$ where $H$ is uniquely determined by the fact that $\Gamma\backslash\Gamma H$ is the closure of the orbit $\Gamma\backslash\Gamma \varphi_g(\SLdR)$. The uniqueness of the limit implies that any subsequence converges, and thus $\rho_r \to \mu_H$ as $r\to 0$.

\section{Examples of cut-and-project sets and their $\SLdR$-closures}\label{sec19}

The first obvious example is when $\scrP$ is a Euclidean lattice $\scrL$ as studied in Sections \ref{sec13}--\ref{sec15}. In this case $m=0$, $G=\SLdR$, $\Gamma=\SLdZ$ and the $\SLdR$ closure of $\scrL$ is the space of lattices $\GamG$; that is $H_g=G$ for any $g$. This observation generalises as follows.

\begin{prop}[\cite{Marklof2013}]\label{prop19.1}
If $m<d$ and $\scrL=\ZZ^{n}g$ is chosen such that $\pi\big|_\scrL$ is injective, then $H_g=G=\SL(n,\RR)$.
\end{prop}

The conditions of this proposition are for instance satisfied in the example studied by Wennberg \cite{Wennberg:2012ev}, where
\begin{equation}
\scrP=\scrQ\times\ZZ\subset\RR^2
\end{equation}
and $\scrQ$ is the one-dimensional cut-and-project set (``Fibonacci quasicrystal'')
\begin{equation}
\scrQ = \bigg\{ \frac{j}{\sqrt{1+\tau^2}} + \frac{1}{\tau\sqrt{1+\tau^2}} \bigg\| \frac{j}{\tau} \bigg\|\bigg\}_{j\in\ZZ}
\end{equation}
where $\tau=\frac{1+\sqrt 5}{2}$ (the golden ratio) and $\|\,\cdot\,\|$ is the distance to the nearest integer. It is an instructive exercise to understand why $\scrQ$ and $\scrP$ are cut-and-project sets as defined in Section \ref{sec18}.

As we shall now see, there are counter examples to the claim in Proposition \ref{prop19.1} for $m\geq d$. 
Probably the most prominent class of {\em quasicrystals} are those constructed from algebraic number fields. The  Penrose tilings fall into this class. Let us briefly sketch how such quasicrystals can be obtained as cut-and-project sets. Let
\begin{itemize}
\item $K$ be a totally real number field of degree $N\geq 2$ over $\QQ$,
\item $\fO_K$ the ring of integers of $K$, and
\item $\pi_1=\id$, $\pi_2,\ldots,\pi_N$ the distinct embeddings $K\hookrightarrow\RR$.
\end{itemize}
We also use $\pi_i$ to denote the component-wise embedddings
\begin{equation}
\begin{split}
\pi_i : & K^d \hookrightarrow \RR^d \\
& \vecx \mapsto (\pi_i(x_1),\ldots ,\pi_i(x_d)) ,
\end{split}
\end{equation}
and similarly for the entry-wise embeddings of $d\times d$ matrices,
\begin{equation}
\pi_i :  \M_d(K) \hookrightarrow \M_d(\RR). 
\end{equation}

Now consider the lattice 
\begin{equation}
\scrL=\{ (\vecx, \pi_2(\vecx),\ldots,\pi_N(\vecx) ) :\vecx\in\fO_K^d \}
\end{equation}
in $\RR^{Nd}$. This is a lattice of full rank. The dimension of the internal space is $m=(N-1)d$. It is a fact of ``basic'' number theory \cite{Weil1974} that $\scrA:=\overline{\pi_\intl(\scrL)} =\RR^m$, so that $\scrV=\RR^{Nd}$. To complete the discussion of this set-up, we need to work out $H_g$ for $g\in G$ and $\delta>0$ so that
\begin{equation}
\scrL=\delta^{1/Nd} \ZZ^{Nd} g.
\end{equation}
(In fact, $\delta=|D_K|^{d/2}$ where $D_K$ is the discriminant of $K$.) The answer is given by the following lemma.

\begin{lem}[\cite{Marklof2013}]
For $g$ as above,
\begin{equation}
H_g = g \SLdR^N g^{-1}
\end{equation}
and
\begin{equation}
\Gamma\cap H_g = g \SL(d,\fO_K) g^{-1} .
\end{equation}
\end{lem}

The group $\SL(d,\fO_K)$ is called the {\em Hilbert modular group}.
The proof of the above lemma is written out in Section 2.2.1 of \cite{Marklof2013}. For a detailed account on how the Penrose tilings fit into this setting, see Section 2.5 of \cite{Marklof2013}.

A further example of a cut-and-project set is to take the union of finite translates of a given cut-and-project set. This is explained in Section 2.3 of \cite{Marklof2013}. Let us here discuss the special case of periodic Delone sets, i.e., the union finite translates of a given lattice $\scrL_0$ of full rank in $\RR^d$. An example of such a set is the honeycomb lattice, which in the context of the Boltzmann-Grad limit of the Lorentz gas was recently studied by Boca et al.~\cite{Boca:2009uz,Boca:2010ty} with different techniques. The scatterer configuration $\scrP$ we are now interested in is the union of $m$ copies of the same lattice $\scrL_0$ translated by $\vect_1,\ldots,\vect_m\in\RR^d$,
\begin{equation}
\scrP = \bigcup_{j=1}^m (\vect_j + \scrL_0) .
\end{equation}
We assume that the $\vect_j$ are chosen in such a way that the above union is disjoint. Let us now show that $\scrP$ can be realised as a cut-and-project set $\scrP(\scrL,\scrW)$. Let
\begin{equation}
\scrL = (\scrL_0 \times \{\vecnull\}) + \sum_{j=1}^m \ZZ\, (\vect_j,\vece_j) \subset\RR^{n},
\end{equation}
where $\vecnull\in\RR^m$ and $\vece_1,\ldots,\vece_m$ are the standard basis vectors in $\RR^m$. The set $\scrL$ is evidently a lattice of full rank in $\RR^{n}$. Note that 
\begin{equation}
\pi_\intl(\scrL) = \sum_{j=1}^m \ZZ\, \vece_j = \ZZ^m,
\end{equation}
and therefore the closure of this set is $\scrA=\ZZ^m$ with connected component $\scrA^0=\{\vecnull\}$.
It follows that for the window set
\begin{equation}
\scrW= \bigcup_{j=1}^m \{ \vece_j \} \subset\scrA
\end{equation}
we indeed have 
\begin{equation}
\scrP(\scrL,\scrW) = \bigcup_{j=1}^m (\vect_j + \scrL_0) .
\end{equation}

Let us now determine $H_g$ in this setting. (Note that the injectivity assumption of Proposition \ref{prop19.1} is now violated.)
Take $g_0\in\SLdR$ so that $\scrL_0=\nbar_0^{-1/d}\ZZ^d g_0$, where $\nbar_0$ is the density of $\scrL_0$. Set
\begin{equation}
T = \begin{pmatrix} \vect_1 \\ \vdots \\ \vect_m \end{pmatrix} \in \M_{m\times d}(\RR).
\end{equation}
We then have
\begin{equation}
\scrL=\nbar_0^{-1/n}\ZZ^n g,
\end{equation}
for
\begin{equation}
g=\nbar_0^{1/n}\begin{pmatrix} \nbar_0^{-1/d} g_0 &  0 \\ T & 1_m \end{pmatrix}
\in\SL(n,\RR).
\end{equation}

If $\scrL_0$ and $\vect_1,\ldots,\vect_m$ are not rationally related, we have the following result for the Ratner subgroup $H_g$. 

\begin{lem}[\cite{Marklof2013}]
Let $\veca_1,\ldots,\veca_d$ be a basis of $\scrL_0$ so that the vectors $\veca_1,\ldots,\veca_d$, $\vect_1,\ldots,\vect_m$ are linearly independent over $\QQ$. Then 
\begin{equation}
H_g = \bigg\{ \begin{pmatrix} h & 0 \\ u & 1_m \end{pmatrix} : h\in\SLdR,\; u\in\M_{m\times d}(\RR) \bigg\} .
\end{equation}
\end{lem}

The Ratner subgroups that appear in the case of rational translates $\vect_j$ are discussed in \cite{Marklof2013}, Section 2.3.1.

\section{Conclusions \label{sec20}}

We have studied several classes of scatterer configurations $\scrP$, for which the dynamics of the Lorentz gas converges in the Boltzmann-Grad limit to a random flight process. This process is described by a {\em generalised} linear Boltzmann equation. Its transition kernel $k(\vecomega',\xi,\vecomega)$ can be obtained from a random point process $\Theta$, which is the limit of a sequence of dilated, randomly rotated copies of $\scrP$ as in \eqref{limitTheta}. If $\scrP$ is sufficiently generic (e.g., a typical realisation of a Poisson point process), then the limit process $\Theta$ is a Poisson point process and the generalised Boltzmann equation reduces to the classical linear Boltzmann equation. If, on the other hand, $\scrP$ is a Euclidean lattice or a quasicrystal, then the limit process $\Theta$ is given by the distribution of random lattices, and the generalised Boltzmann equation does not reduce to the classical case. A striking feature in this case is that the distribution of free path lengths has a power-law tail and no second moment. This is very different to the Poisson case, where the path length distribution is exponential. By taking unions of incommensurate Euclidean lattices, one can in fact obtain path length distributions with any integer power law \cite{Marklof2013b}.

It is an exciting challenge to try to characterise all point processes $\Theta$ that arise as a limit in \eqref{limitTheta} or \eqref{limitKH} for general scatterer configurations $\scrP$. The limit processes we have encountered in these lectures have a common feature: they are invariant under the action of $\SLdR$. That is, the probability of finding a given number of points in each of the test sets $\scrB_1,\ldots,\scrB_k$ is the same as for the test sets $\scrB_1 M ,\ldots,\scrB_k M$, for any $M\in\SLdR$. Although it is not obvious that all limit processes $\Theta$ have this property, I expect that almost all will. A natural objective is therefore:
\begin{equation}
\text{\em Classify all $\SLdR$-invariant point processes in $\RR^d$.}
\end{equation}
This seems an extremely hard problem already in dimension $d=2$, when one looks beyond Ratner's setting of point processes coming from homogeneous spaces (Section \ref{sec12}). Recent breakthroughs include the papers by McMullen \cite{McMullen:2007kn} and Eskin and Mirzakhani \cite{2013arXiv1302.3320E} on the $\SL(2,\RR)$-action on moduli space, which may be mapped to an $\SL(2,\RR)$-invariant point process in $\RR^2$ by analogous arguments as in Section \ref{sec12}, see \cite{Veech:1998ct} for details.  

The origin in $\RR^d$ (which represents the current scatterer location) is a fixed point of the $\SLdR$-action and hence plays a special role. It is natural to focus on those processes $\Theta$ that are independent of the choice of origin. To this end, consider an $\ASL(d,\RR)$-invariant process $\widetilde\Theta$, where $\ASL(d,\RR)$ is the group generated by $\SLdR$ and the group of translations of $\RR^d$. We then obtain the desired $\SLdR$-invariant process $\Theta$ by conditioning $\widetilde\Theta$ to contain the origin. The goal now seems a little easier:
\begin{equation}
\text{\em Classify all $\ASL(d,\RR)$-invariant point processes in $\RR^d$.}
\end{equation}
The limit processes $\Theta$ we have discussed in these lectures fall into this more restricted class. This is evident for the Poisson process, but less obvious in the case of Euclidean lattices and cut-and-project sets. Here the process $\widetilde\Theta$ is constructed via the space of {\em affine} lattices, see \cite{Marklof:2010ib,Marklof2013} for details.

I should point out that the limit processes $\Theta$ in \eqref{limitTheta} or \eqref{limitKH} do not necessarily have to be $\SLdR$-invariant (in the sense defined above), when $\scrP$ is dependent on $r$. Assume for instance that the scatterer locations $\scrP$ are no longer fixed, but oscillate around their equilibrium position at a given points set $\scrP_0$, where the amplitude of oscillation is on the same scale as the scattering radius. We may model this by assuming that the position of the scatterer is a random variable $\vecy + r \vecxi^\vecy$, where $\vecy\in\scrP_0$, and $\{ \vecxi^\vecy: \vecy\in\RR^d\}$ is a random field of identical, independently distributed random vectors $\vecxi^\vecy\in\RR^d$ with rotation-invariant distribution. Denote by $\vecxi^\vecy_\perp\in\RR^{d-1}$ the orthogonal projection onto the plane perpendicular to $\vece_1$. The renormalisation approach in Sections \ref{sec9} and \ref{sec10} shows that the limit random process is given by
\begin{equation}\label{displ}
	\Theta = \big\{ \vecy + \vecxi^\vecy_\perp-\vecxi^\vecnull_\perp : \vecy \in\Theta_0 \big\} ,
\end{equation}
where $\Theta_0$ is the limit process corresponding to the fixed configuration $\scrP_0$. If $\scrP_0$ is given by a Euclidean lattice or cut-and-project set, then $\Theta_0$ is $\SLdR$-invariant but $\Theta$ is not.

\section*{Acknowledgements}

These lecture notes were prepared for the summer schools ``Limit theorems for Dynamical Systems'' at the Bernoulli Center, EPFL Lausanne, 27--31 May 2013, ``Current Topics in Mathematical Physics'' at CIRM Marseille, 1--7 September 2013, and the Masterclass ``Randomness, Classical and Quantum'' at the University of Copenhagen, 4--8 November 2013. I thank the organizers for the invitation to speak at these events, and the participants for stimulating discussions and feedback. Much of the material in these notes is based on joint work with Andreas Str\"ombergsson, and I would like to thank him for the long-standing collaboration. I am grateful to Daniel El-Baz, Jory Griffin, Andreas Str\"ombergsson, Jim Tseng and Ilya Vinogradov for their comments on a first draft of these notes. The research leading to the results presented here has received funding from the European Research Council under the European Union's Seventh Framework Programme (FP/2007-2013) / ERC Grant Agreement n.~291147. I furthermore gratefully acknowledge support from a Royal Society Wolfson Research Merit Award.

\bibliographystyle{amsplain}

\providecommand{\bysame}{\leavevmode\hbox to3em{\hrulefill}\thinspace}
\providecommand{\MR}{\relax\ifhmode\unskip\space\fi MR }
\providecommand{\MRhref}[2]{%
  \href{http://www.ams.org/mathscinet-getitem?mr=#1}{#2}
}
\providecommand{\href}[2]{#2}

\end{document}